\documentclass[11pt]{amsart}
\usepackage{amsfonts}
\usepackage{amssymb}
\newcommand\nota{{\bf  A{{\hskip-6.98pt}\raise1.05pt \hbox{$\boldsymbol /$}\hskip.4pt}}}

\newcommand\notb{{\bf  B{{\hskip-7.25pt}\raise1.05pt \hbox{$\boldsymbol /$}\hskip.4pt} }}

\newcommand\notd{{\bf D  {{\hskip-7.8pt}\raise1.05pt \hbox{$\boldsymbol /$}\hskip1.5pt}}}

\newcommand\noti{{\bf I {{\hskip-5.28pt}\raise1.05pt\hbox{$\boldsymbol /$}} }}

\newcommand\nots{{\bf S {{\hskip-6.98pt}\raise1.05pt \hbox{$\boldsymbol /$}\hskip.8pt} }}

\begin{document}

\newtheorem {theorem}{Theorem}[section]
\newtheorem {lemma}[theorem]{Lemma}
\newtheorem{corollary}[theorem]{Corollary}
\newtheorem {proposition}[theorem]{Proposition}

\theoremstyle{definition}
\newtheorem{example}[theorem]{Example}
\newtheorem{remark}[theorem]{Remark}

\baselineskip18pt

\title{Toward a classification of prime ideals \mbox{in Pr\"ufer domains}}
\author[Marco Fontana]{Marco Fontana$^1$}
\thanks{$^1$Work done under the auspices of GNSAGA (Gruppo Nazionale per le Strutture Algebriche, Geometriche e le loro Applicazioni).}
\address[Marco Fontana]{Dipartimento di Matematica\\ Universit\`a degli Studi Roma Tre\\
    Largo San L. Murialdo, 1\\
    00146 Roma, Italy}
\email{fontana@mat.uniroma3.it}
\author[Evan Houston]{Evan Houston$^2$}
\thanks{$^2$The author was supported by a visiting grant from GNSAGA of INdAM (Istituto Nazionale di Alta Matematica)}
\address[Evan Houston and Thomas G. Lucas]{Department of Mathematics  and Statistics\\
    University of North Carolina at Charlotte\\
    Charlotte, NC 28223 U.S.A.}
\email[Evan Houston]{eghousto@uncc.edu}
\author[Thomas G. Lucas]{Thomas G. Lucas$^3$}
\thanks{$^3$Correspondence to Thomas G. Lucas (tglucas@uncc.edu), 
 Department of Mathematics and Statistics, University of North Carolina Charlotte,  Charlotte, NC 28223, USA}
\email[Thomas G. Lucas]{tglucas@uncc.edu}

\begin{abstract} The primary purpose of this paper is give a classification scheme for the nonzero primes of a Pr\"ufer domain based on five properties. A prime $P$ of a Pr\"ufer domain $R$ could be sharp or not sharp, antesharp or not, divisorial  or not, branched  or unbranched, idempotent  or not. Based on these five basic properties, there  are six types of maximal ideals and twelve types of nonmaximal (nonzero) primes. Both characterizations and examples are given for each type that exists. \end{abstract}

\maketitle
\date{}

\section{Introduction}

Throughout this paper we let  $R$ be an integral  domain with quotient field $K$.  For most of the results, $R$ is assumed to be a Pr\"ufer domain with $R\ne K$.  What we are interested in is a classification scheme for the nonzero primes of $R$. We do this on an individual basis for each prime  rather than considering the collective types a particular Pr\"ufer domain admits in its entire  spectrum of nonzero prime ideals.

Given a nonzero ideal $I$ of  $R$,  we  mainly consider the following algebraic objects associated to $I$:
$$
\begin{array} {rl}
(I:I) :=& \{z\in K \mid zI \subseteq I\},\\
I^{-1} := & (R:I) :=  \{z\in K \mid zI \subseteq R\},\\
I^{v} := & \left(I^{-1}\right)^{-1},\\
T(I) :=& \bigcup \{(R:I^n) \mid n \geq 0\},\\
\Omega(I) :=& \{z\in K \mid \mbox{for each $i\in I$, }  zi^{n(i)} \in R  \mbox{ for some    $n(i)\in \mathbb Z^+$}\} \\
 = & \hskip -8pt ^{(*)} \; \,\bigcap\{R_Q \mid Q \mbox{ ranges over the prime ideals that do not contain } I\},\\
 \Theta(I) :=& \bigcap\{R_N \mid N \mbox{ ranges over the maximal ideals that do not contain } I\},
\end{array}
$$ with $\Omega(I):=K$ [respectively,  $\Theta(I):=K$] if every nonzero prime ideal contains $I$ [respectively,  if every nonzero maximal ideal contains $I$ ].  (For  equation ${(*)}$, see \cite[Theorem 3.2]{FHPcjm}.) For clarity, we use   $T_R(I)$, $\Omega_R(I)$ and $\Theta_R(I)$ when the context involves more than one ring. Obviously, for each nonzero ideal $I$:
$$
R \subseteq (I:I) \subseteq I^{-1} \subseteq T(I) \subseteq \Omega(I) \subseteq \Theta(I).$$
Note that if $P$ is  a nonzero prime of a Pr\"ufer domain $R$, then $(P:P)=R_P\cap \Theta(P)$  (see, for example, \cite[Theorems 3.2 and 3.8]{HuPduals}). In particular, $(M:M)=R$ for each maximal ideal $M$.  By putting together Theorem 3.2, Lemma 3.3 and Proposition 3.9 of \cite{HuPduals}, we have that  for a  nonzero ideal $I$ of a Pr\"ufer domain $R$, $I^{-1}$ a ring implies $I^{-1}=\sqrt I^{-1}=(\sqrt I:\sqrt I)$. This  combination is particularly useful when dealing with primary ideals -- if $Q$ is $P$-primary and $Q^{-1}$ is a ring, then $Q^{-1}=P^{-1}=(P:P)$. 

Recall that a nonzero  ideal $I$ of an integral domain $R$ is said to be divisorial if $I=I^v \, (=(R:(R:I)))$.  It is known that a maximal ideal of a Pr\"ufer domain is divisorial if and only if it is invertible (see, for example,  \cite[Corollary 3.4]{HuPduals} or \cite[Corollary 3.1.3]{FHPbook}).  A characterization for when a nonmaximal (nonzero) prime of a Pr\"ufer domain is divisorial is given in \cite[Proposition 9]{FHPjpaa}. This is  one of the five properties we are interested in. Another is simply whether or not the prime in question is idempotent. A third is whether $P$  is branched, meaning it has a proper $P$-primary ideal, or unbranched (meaning it does not). 
  
In \cite{Gilover}, R. Gilmer introduced the notion of a sharp domain via ``property $(\#)$".  For a domain $R$, he declared $R$ to have  \textit{property $(\#)$} (or to ``satisfy $(\#)$")  if for each pair of  nonempty subsets $\mathcal M$ and $\mathcal N$ of ${\rm Max}(R)$,  $\bigcap \{R_{M_\alpha} \mid M_\alpha\in \mathcal M\} = \bigcap \{R_{N_\beta} \mid N_\beta\in \mathcal N\}$  if and only if $\mathcal M=\mathcal N$.  By default, a local domain is sharp. He proved that an integral domain $R$ (with more than one maximal ideal) is sharp if and only if for each maximal ideal $M_\beta$, $R_{M_\beta}$ does not contain the intersection $\bigcap\{R_{M_\alpha}\mid M_\alpha \in \mathcal M_\beta\}$ where $\mathcal M_\beta={\rm Max}(R)\backslash \{M_\beta\}$ \cite[Lemma 1]{Gilover}. In \cite{LoLu}, K.A.  Loper and the third named author introduced the notion of a maximal ideal $M_\beta$ being \textit{sharp} if $R_{M_\beta}$ does not contain $\bigcap \{R_{M_\alpha} \mid M_\alpha \in \mathcal M_\beta\}$. We extend this definition to say that a nonzero prime $P$ is \textit{sharp} if $R_P$ does not contain $\Theta(P)$.  Whether a (nonzero) prime is sharp or not is the fourth property we include in our classification scheme. 

The fifth property we  consider comes from the characterizations of divisorial nonmaximal primes in Theorem 4.1.10 of \cite{FHPbook},  the proof given for  Theorem 2.5 in \cite{HePrtp} and by several results in \cite{Lurtp}.  We say that a nonzero prime $P$ of an integral domain  $R$ is \textit{antesharp} if each maximal ideal of $(P:P)$ that contains $P$, contracts to $P$ in $R$. If $R$ is a Pr\"ufer domain, each prime of $(P:P)$ is extended from $R$ (see, for example, \cite[Theorem 1]{Gilover} or \cite[Theorem 26.1]{Gilmit}). Hence for $R$ Pr\"ufer,  $P$ is antesharp if and only if it is a maximal ideal of $(P:P)$. In Proposition \ref{pr:1.3}  below we give several other characterizations for antesharp primes in Pr\"ufer domains. Some are rather specific to Pr\"ufer domains, but two involving invertible ideals hold for all integral domains. It is clear that a maximal ideal of a Pr\"ufer domain is antesharp by default. The use of ``antesharp" for this property  comes partially from the fact that if $P$ is a sharp prime of a Pr\"ufer domain, then it is also an antesharp prime  (see Corollary \ref{cor:1.4} below).

Throughout the paper,  we will  draw heavily on Gilmer's characterization of branched and unbranched primes \cite[Theorem 23.3(e)]{Gilmit}. With regard to the sharp condition, we make frequent use of   \cite[Theorem 2]{Gilover} and  \cite[Corollary 2]{GHover}.  Other frequently cited results include  Corollary 3.1.8 and Theorems 3.1.2, 3.3.10 and 4.1.10 of  \cite{FHPbook}. Corollary 3.1.8 is a consequence of Proposition 2.3 and Theorem 3.8 of \cite{HuPduals}.  Theorem 3.1.2 also appears in \cite{HuPduals} as Lemma 3.3 and Theorem 3.2.  Theorem 3.3.10 is from \cite{Hays} and Theorem 4.1.10 is from \cite{FHPjpaa}.   

\section{Sharp and antesharp primes in Pr\"ufer domains}

We begin by giving several ways of characterizing sharp primes and antesharp primes in Pr\"ufer domains. Then we  consider relations between these two properties and the other three.

 \begin{remark}  \label{rk:1.1}  Let $P$ be a nonzero prime of a Pr\"ufer domain $R$. In \cite[Corollary 2]{GHover}, Gilmer and Heinzer proved that (in our terminology), $P$ is sharp if and only if there is a finitely generated ideal $I$ contained in $P$ such that each maximal ideal that contains $I$ also contains $P$ (for $P$ maximal this is \cite[Theorem 2]{GHover}). Also, the only time $P^{-1}$ properly contains $(P:P)$ is when $P$ is invertible \cite[Corollary 3.1.8]{FHPbook}, and therefore maximal. A simple consequence is that in all cases, $(P:P)=R_P\cap \Theta(P)$ \cite[Theorem 3.1.2]{FHPbook}. It follows that  $P$ is sharp if and only if $(P:P)$ is  a proper subset of $\Theta(P)$.  Also, for a nonzero nonmaximal prime ideal $P$,  we know that $P$ is divisorial if and only if $P^{-1} \subsetneq \Theta(P)$ or $(R : \Theta(P)) =P$ \cite[Theorem 4.1.10]{FHPbook}. Hence   if $P$ is a nonzero nonmaximal prime that is sharp, then $P$ is divisorial \cite[Theorem 2.1]{FHPcjm}. 
Finally, note that  $P$ is both branched  and sharp if and only if $P=\sqrt B$  for some finitely generated ideal $B$.  This is simply a combination of the characterization of ``sharp" given  above and Theorem 23.3(e) in \cite{Gilmit} --  specifically, $P$ is branched if and only if  $P$ is minimal over some finitely generated ideal $J \subseteq P$. Hence if  $P$ is both  branched and sharp,  then $P =\sqrt B $ where $B:= I+J$ is such that $P$ is minimal over the finitely generated ideal $J$ and  $I$ is a finitely generated ideal of $R$ such that the only maximal ideals that contain $I$ are those that contain $P$. \end{remark}

  \begin{proposition}  \label{pr:1.2}  Let $R$ be  a Pr\"ufer domain and let  $P$ be a nonzero prime ideal of $R$. 
   Then the following are equivalent.
   \begin{enumerate}
\item [(i)]  $P$ is sharp. 
\item[(ii)] There is a finitely generated ideal $I\subseteq P$ such that the only maximal ideals that contain $I$ are those that contain $P$. 
\item [(iii)] There is a prime ideal $Q\subseteq P$ such that $Q$ is the radical of a finitely generated ideal   with $\Theta(Q)=\Theta(P)$.  
\item [(iv)] There is a prime ideal $Q\subseteq P$ such that $Q$ is the radical of a finitely generated ideal and each maximal ideal that contains $Q$ also contains $P$.
\item [(v)] There is a prime ideal $Q\subseteq P$ such that $Q$ is the radical of a finitely generated ideal and each prime that contains $Q$ is comparable with $P$.
\item [(vi)] There is a prime ideal $Q\subseteq P$ such that $Q$ is the radical of a finitely generated ideal and each ideal that contains $Q$ is comparable with $P$.
\item [(vii)]  There is a finitely generated ideal $I\subseteq P$ such that each ideal that contains $I$ is comparable with $P$.
\end{enumerate}
Moreover, the prime ideal $Q$ in the statements from  {\rm (iii)}  to {\rm  (vi)}   is  sharp and branched.
\end{proposition}  

\begin{proof}  The equivalence of (i) and (ii) is from  \cite[Corollary 2]{GHover}. Clearly (vi) implies (v),  (v) implies (iv),  (iv) implies (iii) and (iii) implies (ii). 

Next we show that (i) implies (vi) and (vii).

Assume $P$ is sharp and let $I\subseteq P$ be a finitely generated ideal such that each maximal ideal that contains $I$ also contains $P$. Let $Q$ be a prime minimal over $I$. Then each maximal ideal $M$ that contains $Q$ also contains $P$.  Since $R_M$ is a valuation domain (or because Spec$(R)$ is treed), $P$ must contain $Q$. It follows that $Q$ is the unique minimal prime of $I$ and therefore, $Q=\sqrt I$. 

For the rest, it suffices to start with an element $r\in R\backslash P$ such that  $rR+P\ne R$, then show that  the ideal $J=rR+I$ contains $P$. For this, we simply see what happens when we localize at a maximal ideal. Clearly, if $N$ is a maximal ideal that does not contain $I$, then $JR_N=R_N=PR_N$. On the other hand, if $M$ is a maximal ideal that contains $I$, then it also contains $P$. Since $r$ is not in $P$ and $PR_M\cap R=P$, $JR_M=IR_M + rR_M=rR_M \supsetneq PR_M$. It follows that  $J \supsetneq P$.  

 To finish we  show (vii) implies  (vi). 
Assume $I$ is a finitely generated ideal of $R$ with $I\subseteq P$ and each ideal  containing $I$ is comparable with $P$. Let $Q$ be minimal over $I$. Since $Q$ is comparable with $P$, we must have $Q\subseteq P$. Since Spec$(R)$ is treed, $Q$ must be the unique minimal prime of $I$ so that $Q=\sqrt I$. 

The last statement is a consequence of Remark \ref{rk:1.1}.  \end{proof}
 We next turn our attention to the antesharp property.

  \begin{proposition}  \label{pr:1.3}  Let $P$ be a  nonzero nonmaximal prime ideal of an integral domain  $R$. Then the following are equivalent. \begin{enumerate}
\item[(i)] $P$ is antesharp (i.e., each maximal ideal of $(P:P)$ that contains $P$ contracts to $P$ in $R$).
\item[(ii)] For each $a\in R\backslash P$, the ideal $A=aR+P$ is invertible.
\item[(iii)] For each prime  $Q$ of $R$ that properly contains $P$,  there is an invertible ideal  $I\subseteq Q$ that properly contains $P$.\end{enumerate}
If $R$ is Pr\"ufer, then (i), (ii) and (iii) are also equivalent to the following.
\begin{enumerate}
\item[(iv)] $P$ is a maximal ideal of $(P:P)$
\item[(v)] Each prime ideal of $(P:P)$ that contains $P$ contracts to $P$ in $R$ and is a maximal ideal of $(P:P)$.
\item[(vi)] For each prime $Q$ of $R$ that properly contains $P$,  there is a finitely generated ideal $I\subseteq Q$ that properly contains $P$. \end{enumerate}
  \end{proposition}

\begin{proof} To see that (i) implies (ii), assume $P$ is antesharp and let $a\in R\backslash P$. Set $A:=aR+P$. As no prime ideal of $R$ that properly contains $P$ survives in $(P:P)$,  we must have  $A(P:P)=(P:P)$. Hence  there are elements $q\in (P:P)$ and $p\in P$ such that $qa+p=1$. It follows that $q=(1-p)/a\in a^{-1}R\bigcap (P:P)$ and thus $q\in A^{-1}$. Since  $p\in A$, $1\in AA^{-1}$ and (i) implies (ii).

Clearly (ii) implies (iii). To finish the general case we show (iii) implies (i). Let $Q$ be a prime ideal that properly contains $P$ and let $I\subseteq Q$ be an invertible ideal that properly contains $P$. Since $P\subsetneq I$, $PI^{-1}\subseteq R$ and $PI^{-1}I=P$ implies $PI^{-1}\subseteq P$.  Hence $I^{-1}\subseteq (P:P)$ and we have $1\in II^{-1}\subseteq I(P:P)\subseteq Q(P:P)$. Thus $Q(P:P)=(P:P)$ and it follows that $P$ is antesharp.

   If $R$ is Pr\"ufer domain, then each prime of $(P:P)$ is uniquely extended from $R$ \cite[Theorem 26.1]{Gilmit}.  Thus (i), (iv) and (v) are equivalent. Also each finitely generated ideal is invertible, so (vi) is equivalent to (iii) finishing  the equivalence of (i)--(vi) in the Pr\"ufer case. 
   \end{proof}
  
If $R$ is not a Pr\"ufer domain, then it is the still the case that (i) implies (vi).  However, the reverse implication does not hold in general, a simple example is the prime $P=(X,Y)$  of  the polynomial ring $K[X,Y,Z]$ where $K$ is  a field.  Also note that since each invertible ideal of a local domain is principal, if $R$ is local, then $P$ is antesharp if and only if it is divided;  i.e., $P$ compares with each principal ideal.

\begin{corollary} \label{cor:1.4}  Let $P$ be a (nonzero) prime ideal of an integral domain   $R$.  \begin{enumerate}
\item[(1)] If $P$ is antesharp and not maximal, then it is divisorial.
\item[(2)] If $R$ is a Pr\"ufer domain and $P$ is sharp (equivalently, $(P:P) \subsetneq \Theta(P)$), then $P$ is a maximal ideal of $(P : P)$ and  antesharp as a prime of $R$.  Moreover, if $P$ is sharp and not maximal, then it is both antesharp and  divisorial.\end{enumerate}
\end{corollary}
  
\begin{proof} Assume $P$ is antesharp and not maximal and let $a\in M\backslash P$ where $M$ is a maximal ideal that contains $P$. Then by statement (ii) of Proposition \ref{pr:1.3}, the ideals $A:=aR+P$ and $B:=a^2R+P$ are both invertible with $B$ properly contained in $A$ since  $A\ne R$. As both $A$ and $B$ are divisorial ideals that properly contain $P$,  $a$ cannot be an element of $P^v$. It follows that $P=P^v$ is divisorial.

If $R$ is a Pr\"ufer domain and $P$ is sharp, then there is a finitely generated ideal $I$ contained in $P$ and in no maximal ideal that does not contain $P$. For each prime $ Q$ properly containing $P$, choose $t \in  Q\setminus P$ and set  $J := tR+I$. For a maximal ideal $M$, if $M$ does not contain $P$, then $PR_M = IR_M = JR_M = R_M$. If $M$ does contain $P$, then we must have $JR_M \supsetneq PR_M$ since $PR_M \cap R = P$ and $t \in J$ does not belong to $ P$. Hence $(Q \supseteq)\ J \supsetneq P$. The conclusions in the first sentence now follow from Proposition \ref{pr:1.3}. Apply statement (1) to verify  the second statement in  (2). 
  \end{proof}
  
     Corollary \ref{cor:1.4} and  Proposition \ref{pr:1.3} suggested the choice of the term ``antesharp" -- a prime $P$ (in a Pr\"ufer domain) is antesharp if it is sharp, and while an antesharp prime $P$ need not contain a finitely generated ideal $I$ contained only in the maximal ideals that contain $P$, each prime that contains $P$ will contain such a finitely generated ideal.

 \bigskip

Several results are already known for prime ideals in Pr\"ufer domains verifying different combinations of the properties of being divisorial, idempotent, branched, sharp and antesharp. Our next goal is to pursue this study to reach  a more complete classification of prime ideals in Pr\"ufer domains. We start with maximal ideals. Recall that if $V$ is a valuation domain with maximal ideal $M$, then either every nonzero ideal is divisorial (equivalently, $M$ is invertible) or $M$ is not divisorial and the other nondivisorial ideals are those of the form $aM$ for some nonzero $a\in V$ (see \cite[Lemma 5.2]{Heinalldiv} and \cite[Lemma 4.2]{BGdplusm}).

\begin{proposition} \label{pr:1.5} Let $R$ be a Pr\"ufer domain and let $M$ be a  branched maximal ideal. Then the following are equivalent.
\begin{enumerate}
\item [(i)] $M$ is sharp but not divisorial.
\item [(ii)]  $M$ is both sharp and  idempotent.
\item [(iii)]  $M$ is idempotent and the radical of a finitely generated ideal.
\item  [(iv)]  $M$ is not divisorial, but divisorial $M$-primary ideals exist.
\item  [(v)]  $M$ is idempotent, but divisorial $M$-primary ideals
exist.
\end{enumerate}
Moreover, if $M$ is idempotent and sharp, then each non-divisorial $M$-primary ideal properly contained in $M$ is of the form $IM$ for some  finitely generated $M$-primary ideal $I$. \end{proposition}

\begin{proof} A (nonzero) maximal ideal of a Pr\"ufer domain is divisorial if and only if it is invertible (for instance, \cite[Lemma 4.1.8]{FHPbook}). Since  $M$ is a branched maximal ideal,  it is minimal over some finitely generated ideal $I$ \cite[Theorem 23.3]{Gilmit}.  Also $M$ is sharp if and only if it contains a finitely generated ideal $J$ that is in no other maximal ideal (Remark \ref{rk:1.1}). In this case,  $M$ is the radical of the finitely generated ideal $I+J$. Also, if $M$ is not idempotent, then necessarily  $MR_M=bR_M$ for some $b\in M$  (\cite[Proposition 5.3.1 (1)]{FHPbook} and \cite[Theorem 17.3 (b)]{Gilmit}). Thus, checking locally,  we have $M=I+bR$, i.e. $M$ is finitely generated and so, in particular, divisorial. Therefore, in our situation, $M$ is either idempotent or divisorial (but not both). It follows that  (i), (ii) and (iii) are equivalent and also that (iv) and (v) are equivalent. Also (iv) implies (i) by \cite[Proposition 2.2]{bruce}.

Assume $M$ is  sharp  but not divisorial. Then $M$ is not finitely generated but there exists a finitely generated (and therefore divisorial) ideal $B$  with $\sqrt B=M$.  Since $M$ is maximal,  $B$ is  $M$-primary. Thus (i) implies both (iv) and (v).

For the ``Moreover" statement, assume $M$ is branched and sharp but not divisorial. By (iii), $M=M^2$ and there is a finitely generated ideal  $B\subsetneq M$  with $\sqrt B=M$.  Let  $Q$ be a proper  $M$-primary ideal that is not divisorial and let $r\in M\backslash Q$. Since $\sqrt Q=M =\sqrt B$ and $B$ is finitely generated, then  $Q$ contains $B^m$ for some positive integer $m$. Since $B$ is $M$-primary, checking locally shows that the finitely generated ideal $I:=B^m+rR$ is a proper $M$-primary ideal that properly contains $Q$.   Since  $Q$ is not divisorial, we may choose $r\in Q^v\backslash Q$ which puts $I$ between $Q$ and $Q^v$. As $I$ is finitely generated and $M$-primary, we must have $Q^v=I$ and  $rR_M =IR_M$. This implies that $QR_M\subsetneq rR_M =IR_M$ and therefore $QR_M\subseteq rMR_M=IMR_M\subsetneq IR_M$ with no ideal properly between $IMR_M$ and $IR_M$ (since $IR_M$ is invertible).  If $Q\ne IM$, there is an $s\in IM\backslash Q$. As $sR_M\subseteq IMR_M \subsetneq  rR_M=IR_M$, the finitely generated ideal $A=B^m+sR$ would be properly between $Q$ and $Q^v=I$ which is absurd. Thus $Q=IM$. \end{proof}

\begin{corollary} \label{cor:1.6} Let $R$ be a Pr\"ufer domain and let $M$ be a  branched  nondivisorial maximal ideal. Then the following are equivalent.
\begin{enumerate}
\item [(i)] $M$ is non-sharp. 
\item  [(ii)]   No $M$-primary ideal is divisorial. 
\item  [(iii)]  $(R:Q)=R$, for each $M$-primary ideal $Q$. 
\end{enumerate}\end{corollary}

\begin{proof}  The equivalence of (i) and (ii) is from the equivalence of statements (i) and (iv) in Proposition \ref{pr:1.5} (see also, \cite[Proposition 2.2]{bruce}).  For the equivalence of (ii) and (iii), assume  $Q$ is  an $M$-primary ideal. Then either $Q^v=R$ (equivalently $(R:Q)=R$) or $Q^v\subseteq M$ is a divisorial  $M$-primary ideal (forcing $M$ to be sharp). 
\end{proof}

Recall from above that for a nonzero ideal $I$, $\Omega(I)=\bigcap \{R_Q\mid Q$ ranges over the prime ideals that do not contain $I\}$ ($=K$, if each nonzero prime contains $I$).  Our next proposition adds to the list of equivalences given in \cite[Theorem 3.3.10]{FHPbook} for having $P^{-1}$ a proper subset of $\Omega(P)$. 

Before stating the proposition we recall three useful facts about duals. First,  for any nonzero ideal $I$ of an integral domain $R$ (Pr\"ufer or not), $(II^{-1})^{-1}$ is always a ring and always equal to $(II^{-1}:II^{-1})$ \cite[Proposition 7.2]{Bass}. For the other two, start with a nonzero prime $P$ of a Pr\"ufer domain such that $P$ is the radical of a finitely generated ideal and an ideal $I$ such that $P$ is minimal over $I$ with $IR_P\subsetneq PR_P$. Then $I^{-1}$ is not a ring and each prime that properly contains $P$ blows up in $P^{-1}$ \cite[Lemmas 21 and 22]{Lurtp} (see also \cite[Lemmas 4.2.25 and 4.2.26]{FHPbook}).

\begin{proposition} \label{pr:1.7} Let $P$ be a nonzero branched prime of a Pr\"ufer domain
$R$. Then the following are equivalent.
\begin{enumerate} 
\item[(i)] $P$ is sharp. 
\item [(ii)] $P$ is  the radical of a finitely generated ideal.
\item  [(iii)] $P\Omega(P)=\Omega(P)$.
\item   [(iv)] $P^{-1}\subsetneq \Omega(P)$. 
\item   [(v)] If $Q$ is a  proper $P$-primary ideal, then\\
 $\mbox{   \: \:   } QQ^{-1}=P$, whenever $P$
is not maximal  \  and\\ $\mbox{   \; \;   }  QQ^{-1}\supseteq P$,  whenever $P$ is maximal. 
\item  [(vi)] There exists a proper $P$-primary ideal that is divisorial.
\end{enumerate}\end{proposition}

\begin{proof}   The equivalence of (i) and (ii) is from Remark \ref{rk:1.1}, and the equivalence of (ii), (iii) and (iv)  is part of \cite[Theorem 3.3.10]{FHPbook}. 

Of the remaining implications, the one that is easiest to establish is  (vi)
implies (iv). Since $\Omega(P)$ is a ring that contains $P^{-1}$, there is
nothing to prove if $P$ is invertible, because in this case $P^{-1}$ is not a ring \cite[Theorem 3.1.2 and Corollary 3.1.3]{FHPbook}. Thus we may assume $P^{-1}=(P:P)$ \cite[Corollary 3.1.8]{FHPbook} and
that there is a proper $P$-primary ideal $Q$ such that $Q=Q^v$.   If $N$ is a prime that does not contain $P$, then $Q^{-1}$ is contained in $R_N$. Thus $\Omega(P)$ contains $Q^{-1}$. It follows that $P^{-1}\subsetneq  Q^{-1}\subseteq \Omega(P)$.

It is nearly as easy to show that (v) implies (iv).  As with the proof of (vi)
implies (iv),  if some proper $P$-primary ideal $Q$ is invertible, then
$P^{-1}\subsetneq Q^{-1}\subsetneq \Omega(P)$ since $\Omega(P)$ is a ring that contains $Q^{-1}$ (which is not a ring, since $Q$ is invertible).  Hence to show  (v) implies (iv)  we may assume each proper $P$-primary ideal $Q$ is such that $QQ^{-1}=P$. Since we have also  assumed $P$ is branched, there are proper $P$-primary ideals. Let $Q$ be one such ideal. If $P^{-1}=\Omega(P)$, then $Q^{-1}=\Omega(P)$ as well since it is always the case that  $P^{-1}\subseteq Q^{-1}\subseteq \Omega(P)$ (no matter whether $R$ is a Pr\"ufer or not).
Since $P^{-1}$ is a ring, we have  $P^{-1} = R_P \cap \Theta(P)$ \cite[Theorem 3.1.2]{FHPbook}. Therefore, $P= QQ^{-1} = QP^{-1} \subseteq QR_P \cap R =Q$, leading to the contradiction that $Q=P$. Thus we have that $P^{-1}$ is properly contained in $\Omega(P)$.

We next show that (ii) implies (v).  Assume $P=\sqrt I$ where $I$ is a finitely
generated ideal of $R$.  Then $P\Omega(P)=\Omega(P)$ (since we know already that (ii)$\Leftrightarrow$(iii)) and therefore, $(P:P)$ is a proper subring of $\Omega(P)$.  Since $R$ is a Pr\"ufer domain, each prime of $(P:P)$  is extended from $R$. Also, by \cite[Lemma 4.2.26]{FHPbook}, each prime that properly contains $P$ blows up in $(P:P)$. Thus $P$ is a maximal ideal of $(P:P)$. 

Let $Q$ be a proper $P$-primary ideal.  Then $Q^{-1}$ is contained in
$\Omega(P)$, and this containment must be proper since $Q^{-1}$ is a ring
if and only if it equals $(P:P)$ \cite[Proposition 3.1.16]{FHPbook},  already known to be  a proper subring of $\Omega(P)$.

We first consider the case that $P$ is maximal. If $P$ is invertible, then
so is $Q$ \cite[Theorem 23.3(b)]{Gilmit}. Thus we may assume $P$ is not invertible.  In this case $P^{-1}=(P:P)=R$ \cite[Corollary 3.1.8]{FHPbook}.  No maximal ideal other than $P$ can contain $QQ^{-1}$
so either $Q$ is invertible or $QQ^{-1}$ is contained in $P$, making it a
$P$-primary ideal whose dual is a ring \cite[Proposition 3.1.1 (2)]{FHPbook}. In the latter case, apply  \cite[Lemma 4.2.25]{FHPbook} to see that we must have $QQ^{-1}=P$.

Now assume $P$ is not maximal. In this case, $P^{-1}=(P:P)$ with $P$ a
maximal ideal of $(P:P)$ \cite[Corollary 3.1.8 and Lemma 4.2.26]{FHPbook}.  It is easy to show as above that $QP^{-1}\subseteq Q$. Therefore, in particular,  $QQ^{-1}P^{-1}\subseteq QQ^{-1}$. Thus   $P^{-1} \subseteq (R:QQ^{-1})$ and we know that $(R:QQ^{-1}) \subseteq Q^{-1}\subsetneq \Omega(P)$. Since $P$ is divisorial \cite[Corollary 4.1.12]{FHPbook}, we have $Q\subseteq QQ^{-1}\subseteq  (QQ^{-1})^v \subseteq P$.  Also, since there are no rings properly between $(P:P)=P^{-1}$ and $\Omega(P)$ \cite[Theorem 3.3.7]{FHPbook},  we have $P^{-1} =(QQ^{-1}:QQ^{-1}) =(R:QQ^{-1})$. As $P$ is a maximal ideal of $P^{-1}$, $QQ^{-1}$ is a $P$-primary ideal in $P^{-1}$. By applying  \cite[Lemma 4.2.25]{FHPbook},
we deduce that  $QQ^{-1}=P$.

To finish we show that (ii) implies (vi). Let $B$ be a finitely generated ideal of $R$ such that $P=\sqrt B$.  If $P$ is maximal, then $B^2$ is proper $P$-primary ideal that is divisorial. On the other hand, if $P$ is not maximal, then we at least have that it is a maximal ideal of $P^{-1}=(P:P)$, since  no prime that properly contains $P$ survives in $P^{-1}$ \cite[Lemma 4.2.26]{FHPbook}.  In this case, we consider the ideal $Q:=B^2P^{-1}$. Since $P^{-1}$ is a ring and $P$ is a maximal ideal of $P^{-1}$, $Q\ne P$ and, clearly, it is $P$-primary in  $P^{-1}$ and, therefore, in  $R$ (because $Q \subsetneq PP^{-1}=P \subsetneq R$).  Since $B$ is invertible as an ideal of $R$ and $P^{-1}$ is a divisorial fractional ideal of $R$, $Q$ is also a divisorial ideal of $R$. 
\end{proof}

Even more is true for nonmaximal branched primes. The following rather odd lemma is the key for a short proof that if $P$ is a nonmaximal branched prime, then it is sharp if and only if each $P$-primary ideal is divisorial.

\begin{lemma} \label {le:1.8}
 Let $R$ be a Pr\"ufer domain and let $P$ be a nonzero nonmaximal branched prime of $R$ that is sharp. If  $Q$ is a proper $P$-primary ideal of $R$ and $B$ is an ideal of $R$ such that $Q \subsetneq B \subseteq P$, then there is a finitely generated ideal $I$ of $R$ such that $Q \subsetneq I \subsetneq B$. 
\end{lemma}

\begin{proof} Let $Q$ be a proper $P$-primary ideal and let $B$ be such that $Q\subsetneq B\subseteq P$.  Since $P$ is sharp and branched, there is a finitely generated ideal $J$ such that $\sqrt J=P$ (Remark \ref{rk:1.1}). Since $\sqrt Q=P$ and $J$ is finitely generated,  $J^m\subsetneq Q$ for  some $m \geq 1$. Choose an element $t\in B\backslash Q$, then  set $A:=tR+J^m$. Let $N$ be a maximal ideal of $R$. If $N$ does not contain $P$, then  $R_N=J^mR_N=QR_N=AR_N$. On the other hand, if $N$ contains $P$, then $QR_N\cap R=Q$ since $Q$ is $P$-primary and therefore $AR_N$  properly contains $QR_N$, since $t \in AR_N  \setminus  QR_N$ and $R_N$ is a valuation domain. Hence $ Q \subsetneq A \subseteq B$. 

Let $M$ be a maximal ideal that contains $P$ and let $s\in M\backslash P$.  Then $ts\in P\backslash Q$ since $Q$ is $P$-primary, $t\in P\backslash Q$ and $s\in M\backslash P$.  Thus $tsR_M$ properly contains $(QR_M \supseteq)\ Q$. Let $I:=tsR+J^m$, obviously $I \subseteq A$. Since $J^m \subsetneq Q$, $ts \in P\setminus Q$ and $s \in M \setminus P$
we have $QR_M \subsetneq tsR_M= tsR_M+ J^mR_M = IR_M \subsetneq  tR_M= tR_M + J^mR_M =AR_M$.  Let $N$ be a maximal ideal distinct from $M$. If  $N$  does not contain $P$, we have $R_N=J^mR_N=QR_N=AR_N=IR_N$.  On the other hand, if $N$ does contain $P$, we  have that $tsR_N$ properly contains $QR_N$ since $QR_N\cap R=Q$ and $ts\notin Q$. It follows that  $ Q \subsetneq I \ (\subseteq A)$. Since, as we remarked above,  $IR_M$ is a proper subset of $AR_M$, we conclude that $ Q \subsetneq I \subsetneq A \subseteq B$. \end{proof}

Using Lemma \ref{le:1.8} we have the following extension of Proposition \ref{pr:1.7} to nonmaximal primes.

\begin{proposition} \label {pr:1.9} Let $P$ be a nonzero prime of a Pr\"ufer domain $R$. If $P$ is not maximal, then the following are equivalent. \begin{enumerate}
\item[(i)] $P$ is both branched and sharp.
\item[(ii)] $P$  is  the radical of a finitely generated ideal.
\item[(iii)] There exists a proper  $P$-primary ideal that is divisorial.
\item[(iv)] Proper $P$-primary ideals exist and each (proper) $P$-primary ideal is divisorial. \end{enumerate}\end{proposition}

\begin{proof} Since $P$ is branched (i.e., proper $P$-primary ideals exist) if and only if it is minimal over a finitely generated ideal, the first three statements are equivalent by  Proposition \ref{pr:1.7}. Clearly (iv) implies (iii). Thus to complete the proof it suffices to show that (i) implies (iv). 

Assume $P$ is a nonzero nonmaximal  prime ideal that is both branched and sharp. Then $P$ is divisorial by  Remark \ref{rk:1.1}.  Let  $Q$ be a proper $P$-primary ideal.  Since $P$ is divisorial, $Q\subseteq Q^v\subseteq P$. By Lemma \ref{le:1.8}, if there is an element $t\in Q^v\backslash Q$, then there is a finitely generated (= invertible ideal) $I$ such that  $Q\subsetneq I=I^v\subsetneq  Q^v$ which is impossible. Hence $Q=Q^v$. \end{proof}

 For the final results of this section we consider unbranched maximal ideals.

\begin{proposition}  \label{pr:1.10} Let $R$ be a Pr\"ufer domain and let $M$ be an  unbranched maximal ideal. Then the following are equivalent:
\begin{enumerate}
\item [(i)] $M$ is sharp.
\item  [(ii)]  There is a  divisorial prime $P$ contained in $M$ that is contained in no  other maximal ideal.  
\item [(iii)]  There is a  divisorial ideal contained in $M$ that is contained in no  maximal ideal.
\item [(iv)]  There is a sharp prime contained in $M$ that is contained in no other maximal ideal. 
\item [(v)]  There is a sharp branched prime contained in $M$ that is contained in no other maximal ideal. 
\item[(vi)] There is a chain of  sharp branched primes $\{P_\alpha\}$ such that $\bigcup P_\alpha=M$ and $R/P_\alpha$ is a valuation domain for each $P_\alpha$.\end{enumerate}\end{proposition}

\begin{proof}  Clearly, (v) implies (iv). The equivalence of (i) and (iii) is by \cite[Proposition 2.2]{bruce}, as is the fact that (ii) implies (i). 

If $M$ is sharp, then there is a nonzero  finitely generated ideal $J\subsetneq M$ such that $M$ is the only maximal ideal containing $J$ \cite[Theorem 2]{Gilover}. Let $P$ be a  prime minimal over $J$. Then $P$ is branched by  \cite[Theorem 23.3(e)]{Gilmit}. Also  $M$ is the only maximal ideal that contains $P$ and $P=\sqrt J$. Thus $P$ is sharp, branched and not maximal. Therefore $P$ is divisorial by Corollary \ref{cor:1.4}.  Hence (i) implies both (ii) and  (v). 

Continuing with the notation of the previous paragraph, $R/P$ must be  a valuation domain with $M/P$ unbranched. Thus  $M=\bigcup Q_\beta$ where the $Q_\beta$'s range over the nonmaximal primes of $R$ that contain $P$. Each $Q_\beta$ is contained in $M$ and no other maximal ideal. Also, each  $Q_\beta$ is sharp by Proposition \ref{pr:1.2}.   The subfamily $\{P_\alpha\}$ consisting of the branched members of $\{Q_\beta\}$ is such that $\bigcup P_\alpha=M$. Moreover, each $R/P_\alpha$ is a valuation domain. Thus (i) implies (vi).

To see that (vi) implies (v), simply note that if $\{P_\alpha\}$ is a chain of  sharp branched primes such that $\bigcup P_\alpha =M$ with $R/P_\alpha$ a valuation domain for each $P_\alpha$, then any one of the $P_\alpha$'s will satisfy the statement in (v). 

Finally assume there is a prime $P\subseteq M$ that is sharp. Then there is a finitely generated ideal $I\subseteq P$ such that each maximal ideal that contains $I$ also contains $P$. If $M$ is the only maximal ideal containing $P$, it is also the only maximal ideal containing $I$. So (iv) implies (i) by Remark \ref{rk:1.1}.\end{proof}

The Pr\"ufer domain in \cite[Example 38]{Lurtp} has an unbranched maximal ideal that is a union of sharp branched primes, but is not itself sharp.

\section{Classifying prime ideals}

For a nonzero prime $P$ of an integral domain $R$, we let $\boldsymbol{\Lambda}(P)=\langle{\bf V,W,X,Y,Z}\rangle$ with the five positions  
 $ \langle\mbox{\bf \mbox{--}, \mbox{--},  \mbox{--},  \mbox{--},  \mbox{--}}   \rangle$  \ corresponding to whether the prime is, respectively:  sharp (for short ``${\bf S}$''), antesharp  (for short ``${\bf A}$''),  divisorial (for short ``{\bf D}''),  branched (for short ``${\bf B}$''), idempotent (for short  ``${\bf I}$'') or not (for short, respectively, ``$\nots$, $\nota$, $\notd$,  $\notb$, $\noti$''). In  other terms,  ${\bf V} \in \{{\bf S}, \nots\}$,  ${\bf W}\in \{{\bf A}, \nota\}$, ${\bf X} \in \{{\bf D}, \notd\}$,  ${\bf Y} \in  \{{\bf B}, \notb\}$ and ${\bf Z} \in  \{{\bf I}, \noti\}$. Since each maximal ideal is antesharp in a Pr\"ufer domain, we drop the reference to antesharp and use only the 4-tuple $\boldsymbol{\Lambda}(M)=\langle{\bf V,X,Y,Z}\rangle$ for the maximal ideal $M$.  (The notation $\boldsymbol{\Lambda}_R(P)$ will be used when the context involves more than one ring.)

We first consider the possible values of $\boldsymbol{\Lambda}(M)$ when $M$ is a maximal ideal (of a Pr\"ufer domain). Of the sixteen  4-tuples, exactly six can occur: 

$$\begin{aligned}
\langle {\bf S},{\bf D},{\bf B},\noti\rangle  & \quad &  
\langle {\bf S},\notd,{\bf B},{\bf I}\rangle   &  \quad & \langle {\bf S},\notd,\notb,{\bf I}\rangle& \, \phantom{.} \\
\langle \nots,\notd,{\bf B},\noti\rangle \hskip2pt  &  \quad & \langle \nots,\notd,{\bf B},{\bf I}\rangle & \quad   &  \langle \nots,\notd,\notb,{\bf I}\rangle & \, .
\end{aligned}$$

For now we only characterize each of these six types and show that there can be no others. In Section \ref{Examples} we give examples of each type.

\begin{theorem} \label{th:1.11} Let $M$ be a nonzero maximal ideal of a Pr\"ufer domain $R$. Then the only possible values for $\boldsymbol{\Lambda}(M)$ are 
$\langle {\bf S},{\bf D},{\bf B},\noti\rangle$,  
$\langle {\bf S},\notd,{\bf B},{\bf I}\rangle$, 
 $\langle {\bf S},\notd,\notb,{\bf I}\rangle$, 
$\langle \nots,\notd,{\bf B},{\bf I}\rangle$, 
$\langle \nots,\notd,{\bf B},\noti\rangle$  and 
$\langle \nots,\notd,\notb,{\bf I}\rangle$. 

More precisely,

\begin{enumerate}

\item[(1)] $\boldsymbol{\Lambda}(M)=\langle {\bf S},{\bf D},{\bf B},\noti\rangle$ if and only if  $R \subsetneq M^{-1} \; \mbox{ (i.e. $M$ is invertible)}$.

\item[(2)] The following are equivalent:\\
 {\rm (2-i)}  $\boldsymbol{\Lambda}(M)=\langle {\bf S},\notd,{\bf B},{\bf I}\rangle$. \newline
{\rm (2-ii)}   $M$ is not finitely generated, but there is a finitely generated ideal $I$ such that $M=\sqrt I$.\newline
{\rm (2-iii)}   $M$ is not locally principal, but $M=\sqrt I$ for some finitely generated $I$. \newline
{\rm (2-iv)}   $M$ is not invertible, but there are proper $M$-primary ideals that are divisorial.\\
In this case, $R = M^{-1}=T(M)\subsetneq \Omega(M)\subseteq  \Theta(M)$.

\item[(3)]  The following are equivalent:\\
{\rm (3-i)}  $\boldsymbol{\Lambda}(M)=\langle {\bf S},\notd,\notb,{\bf I}\rangle$.
\newline
 {\rm (3-ii)}   $M=\bigcup \{P_\alpha \mid  P_\alpha \subsetneq M,\; P_\alpha \in \mbox{\rm Spec}(R)\}$ and $R\subsetneq \Theta(M)$.
\newline
 {\rm (3-iii)}    $M$ has no proper primary ideals and there is a divisorial ideal contained in $M$ and no other maximal ideal.
\newline
 {\rm (3-iv)}    $M$ has no proper primary ideals and there is a finitely generated ideal contained in $M$ and no other maximal ideal. 
\newline
 {\rm (3-v)}   $M$ has no proper primary ideals and there is a sharp prime  contained in $M$ and no other maximal ideal.\newline
  In this case, $R = M^{-1}=T(M)= \Omega(M) \subsetneq \Theta(M)$.

\item[(4)]  The following are equivalent:\\
 {\rm (4-i)}  $\boldsymbol{\Lambda}(M)=\langle \nots,\notd,{\bf B},{\bf I}\rangle$.  \newline
{\rm (4-ii)}  $M$ is minimal over a finitely generated ideal, but neither the radical of a finitely generated ideal nor locally principal.
\newline
{\rm (4-iii)}  there are proper $M$-primary ideals, none of these are divisorial and $M$ is not locally principal. \\
{\rm (4-iv)}  $M$ is the radical of an ideal properly contained in $M$, no $M$-primary ideal is divisorial and $M$ is not locally principal. \\
In this case, $R = M^{-1}=T(M)= \Omega(M)=\Theta(M)$.

 \item[(5)]  The following are equivalent:\\
 {\rm (5-i)}  $\boldsymbol{\Lambda}(M)=\langle \nots,\notd,{\bf B},\noti\rangle$. 
\newline
 {\rm (5-ii)}    $M$ is locally principal, but it is not the radical of a finitely generated ideal. 
\newline
 {\rm (5-iii)}   $M$ is locally principal, but no proper $M$-primary ideal is divisorial. \\
{\rm (5-iv)}  $M$ is  the radical of an ideal properly contained in $M$, no $M$-primary ideal is divisorial and $M$ is locally principal. \\ 
 In this case, $R = M^{-1}=T(M)= \Omega(M)=\Theta(M)$.

\item[(6)]  The following are equivalent:\\
{\rm (6-i)}  $\boldsymbol{\Lambda}(M)=\langle \nots,\notd,\notb,{\bf I}\rangle$. 
\newline
{\rm (6-ii)}  $M=\bigcup \{P_\alpha \mid P_\alpha \subsetneq M ,\; P_\alpha \in \mbox{\rm Spec}(R)\}$  and $R=\Theta (M)$.
\newline
{\rm (6-iii)}   $M$ has no proper primary ideals and each finitely generated ideal contained in $M$ is contained in at least one other  maximal ideal.\\
{\rm (6-iv)}  $M$ has no proper primary ideals and each divisorial ideal contained in $M$ is contained in at least one other maximal ideal.
\newline
  In this case, $R = M^{-1}=T(M)= \Omega(M) = \Theta(M)$.  
\end{enumerate}
\end{theorem} 

\begin{proof} Since $R$ is Pr\"ufer, $M$ is divisorial if and only if it is invertible.  In such a case $M$ is branched, sharp and not idempotent.
Thus $\langle {\bf S},{\bf D}, {\bf B},\noti\rangle$ is the only possiblity with ``${\bf D}$". 
Also note that if $M$ is unbranched, then it must be idempotent. So $\boldsymbol{\Lambda}(M)$ cannot contain both $\notb$ and $\noti$. Moreover, by Proposition \ref{pr:1.5}, $\boldsymbol{\Lambda} (M)=\langle {\bf S},\notd, {\bf B},\noti\rangle$ is also impossible. All that are left are the six listed above.

(1) We have observed that $\boldsymbol{\Lambda}(M) =
 \langle {\bf S},{\bf D},{\bf B},\noti\rangle$ holds if and only if $M$ is invertible and this happens if and only if
 $R \subsetneq M^{-1}$  (for instance, \cite[Corollary 3.1.3]{FHPbook}).
 
 (2) The equivalence of (2-i) and (2-ii) follows easily from 
 Remark \ref{rk:1.1}, Proposition \ref{pr:1.5}  and the fact that  ``not divisorial'' and ``not finitely generated'' are equivalent notions for maximal ideals in Pr\"ufer domains.   Clearly (2-ii)  and (2-iv) are equivalent from Proposition \ref{pr:1.5}.  It is obvious that   (2-iii) implies (2-ii). Finally (2-i) implies (2-iii) by Proposition \ref{pr:1.5} and the fact that in a Pr\"ufer domain, a maximal ideal is idempotent if and only if it is not locally principal.

Since $M=M^2$, it is trivial that $M^{-1}=T(M)$. The rest of the last statement follows from \cite[Theorem 3.3.10]{FHPbook} and the assumption that $M$ is sharp.

(3)  Since $R$ is a Pr\"ufer domain, $M$ is unbranched if and only if it is the union of the primes properly contained in $M$ \cite[Theorem 23.3]{Gilmit}. Thus (3-i) and (3-ii) are equivalent. The equivalence of  (3-v) and  (3-iv) is from Proposition \ref{pr:1.2}. Also  from Proposition \ref{pr:1.2}, we have that (3-i) implies (3-v). It is obvious that  (3-iv) implies (3-iii).  To finish the equivalence, we have  that (3-iii) implies (3-i) by  Proposition \ref{pr:1.10}. 

As in (2), the last statement follows from \cite[Theorem 3.3.10]{FHPbook} and from the fact that $M^{-1}=T(M)$ when $M=M^2$.

 (4)  Since $R$ is a Pr\"ufer domain,  $M$  is locally principal if and only if   $M\neq M^2$.  That (4-i) implies (4-ii), (4-iii) and (4-iv) now follows from  Corollary \ref{cor:1.6}  and \cite[Theorem 23.3]{Gilmit}.  Since $M$ is maximal, the only real difference between the statements in  (4-iii)  and (4-iv) is that in (4-iii) only the proper $M$-primary are required to be nondivisorial but $M$ not locally principal implies it is not invertible, so it is not divisorial. To see that (4-iii) implies (4-ii) and (4-i), again apply Theorem 23.3 of \cite{Gilmit} to see that $M$ is minimal over a finitely generated ideal. Next, use Corollary \ref{cor:1.6} and Remark \ref{rk:1.1} to conclude that $M$ is  neither sharp nor the radical of a finitely generated ideal.  Finally  to see that (4-ii) implies (4-i), the particular set of restrictions in (4-ii) clearly implies $M$ is branched, idempotent and not divisorial. Thus by (2), it is sharp if and only if it is the radical of a finitely generated ideal. Hence (4-ii) implies (4-i). 

The last statement always holds  in a nonsharp case. So it holds here and in  (5) and (6) as well.

  (5)   By  Corollary \ref{cor:1.6}, \cite[Theorem 23.3(e)]{Gilmit} and (3) we have that  (5-i)  implies (5-iv).   It is obvious that  (5-iv) implies (5-iii) and 
 (5-iii) implies (5-ii). (5-ii) implies (5-i)  by Corollary \ref{cor:1.6}, since $M$ is neither idempotent nor divisorial.

(6) As in (3), we use the fact that there are no proper $M$-primary ideals if  and only if $M$  is the union of the primes that are properly contained in $M$ \cite[Theorem 23.3]{Gilmit}. Also, $M$ fails to be sharp if and only if each finitely generated ideal it contains is contained in at least one other maximal ideal \cite[Theorem 2]{GHover}. Since an unbranched prime must be idempotent, the equivalence of the first three statements is now clear. Clearly, (6-iv) implies (6-iii). Finally,  (6-i) implies (6-iv) by  Proposition \ref{pr:1.10}.\end{proof}

For nonmaximal prime ideals $P$ of a Pr\"ufer domain, the possible values for $ \boldsymbol{\Lambda}(P)$ are described next. In Section \ref{Examples} we will show with explicit constructions that all the possible values  for $ \boldsymbol{\Lambda}(P)$  and for $ \boldsymbol{\Lambda}(M)$ are effectively attained. Of the thirty-two  possible 5-tuples, exactly twelve can (and do) occur:

  $$\begin{aligned} 
  \langle {\bf S},{\bf A}, {\bf D},{\bf B},{\bf I}\rangle  & \quad &  
 \langle {\bf S},{\bf A},{\bf D},{\bf B},\noti\rangle   & \quad & 
 \langle {\bf S},{\bf A},{\bf D},\notb,{\bf I}\rangle  & \, \phantom{.}
\\  
  \langle \nots,{\bf A},{\bf D},{\bf B},{\bf I}\rangle    & \quad &   
 \langle  \nots,{\bf A},{\bf D},{\bf B},\noti\rangle   & \quad & 
 \langle  \nots,{\bf A},{\bf D},\notb,{\bf I}\rangle & \, \phantom{.}
\\  
  \langle \nots,\nota,{\bf D},{\bf B},{\bf I}\rangle   & \quad &   
 \langle \nots,\nota,{\bf D},{\bf B},\noti \rangle   & \quad &  
 \langle \nots,\nota,{\bf D},\notb,{\bf I}\rangle & \, \phantom{.}
\\  
  \langle \nots,\nota,\notd,{\bf B},{\bf I}\rangle    & \quad &  
 \langle \nots,\nota,\notd,{\bf B},\noti\rangle    & \quad &  
 \langle \nots,\nota,\notd,\notb,{\bf I} \rangle & \, . \end{aligned}$$

We say that a nonzero prime $P$ verifies 
$(\boldsymbol{\natural})$  if for each prime $Q$ properly containing $P$, there is a finitely generated ideal $I$ of $R$ such that $P \subsetneq I \subseteq Q$. By Proposition \ref{pr:1.3} this is equivalent to saying that $P$ is antesharp.

\begin{theorem} \label{th:1.12} Let $P$ be a nonzero nonmaximal prime  ideal of a Pr\"ufer domain. Then   the only possible values for \ $ \boldsymbol{\Lambda}(P)$ \ are\;
$\langle {\bf S},{\bf A},{\bf D},{\bf B},{\bf I}\rangle$, 
$\langle {\bf S},{\bf A},{\bf D},{\bf B},\noti\rangle$,  
$\langle {\bf S},{\bf A},{\bf D},\notb, {\bf I}\rangle$,\,
$\langle \nots,{\bf A},{\bf D},{\bf B},{\bf I}\rangle$, \,
$\langle \nots,{\bf A},{\bf D},{\bf B},\noti\rangle$, \,
$\langle \nots,{\bf A},{\bf D},\notb,{\bf I}\rangle$, \,
$\langle \nots,\nota,{\bf D},{\bf B}, {\bf I}\rangle$, \,
$\langle \nots,\nota,{\bf D},{\bf B},\noti\rangle$, 
$\langle \nots,\nota,{\bf D},\notb, {\bf I}\rangle$, 
$\langle \nots,\nota,\notd,{\bf B},{\bf I} \rangle$,  
$\langle \nots,\nota,\notd,{\bf B}, \noti\rangle$, and $\langle \nots,\nota,\notd,\notb,{\bf I}\rangle$.
More precisely: 
\begin{enumerate}

\item The following are equivalent: \\
{\rm (1-i)}  $\boldsymbol{\Lambda} (P) =\langle {\bf S},{\bf A},{\bf D},{\bf B},{\bf I}\rangle$. \newline
{\rm (1-ii)}  $PR_P$ is not principal and $P$ is the radical of a finitely generated ideal. \newline 
{\rm (1-iii)}  $PR_P$ is not principal and some proper $P$-primary ideal is divisorial. \\
{\rm (1-iv)}  $PR_P$ is not principal, proper $P$-primary ideals exist and each $P$-primary ideal is divisorial.

\item[(2)] The following are equivalent: \\
  {\rm (2-i)} 
$\boldsymbol{\Lambda} (P) = \langle {\bf S},{\bf A},{\bf D},{\bf B},\noti\rangle$. \newline
{\rm (2-ii)}   $PR_P$ is  principal and $P$ is the radical of a finitely generated ideal. \newline 
{\rm (2-iii)}  $PR_P$ is  principal and some proper $P$-primary ideal is divisorial. \newline
{\rm (2-iv)}  $PR_P$ is  principal and each $P$-primary ideal is divisorial. \newline
{\rm (2-v)} $P$ is an invertible ideal of $(P:P)$.

\item[(3)] The following are equivalent:\\
{\rm (3-i)} 
$\boldsymbol{\Lambda} (P) =\langle {\bf S},{\bf A},{\bf D},\notb, {\bf I}\rangle$.
 \newline  
{\rm (3-ii)}  $P$ has no proper $P$-primary ideals and there is a finitely generated ideal $I\subseteq P$ such that each ideal that contains $I$ is comparable with $P$. 
\newline 
{\rm (3-iii)}   $P$ has no proper $P$-primary ideals and there is a prime $Q\subsetneq P$ such that $Q$ is the radical of a finitely generated ideal and each maximal (prime) ideal that contains $Q$, also contains (is comparable with) $P$.

 \item[(4)] The following are equivalent: \\
{\rm (4-i)} 
$\boldsymbol{\Lambda} (P) =\langle \nots,{\bf A},{\bf D},{\bf B}, {\bf I}\rangle$. \newline
{\rm (4-ii)}  $PR_P$ is not  principal, proper $P$-primary ideals exist, but $P$ is the only divisorial $P$-primary ideal and $P$ verifies $(\boldsymbol{\natural})$. \newline 
{\rm (4-iii)}  $PR_P$ is not  principal, $P$ is  minimal over some finitely generated ideal, but each finitely generated ideal contained in $P$ is contained in at least one maximal ideal that does not contain $P$ and $P$ verifies $(\boldsymbol{\natural})$. \newline 
{\rm (4-iv)}  $PR_P$ is not principal, $P$ is  minimal over some finitely generated ideal, but each finitely generated ideal contained in $P$ is contained infinitely many maximal ideals that do not contain $P$ and $P$ verifies $(\boldsymbol{\natural})$.

\item[(5)] The following are equivalent: \\
{\rm (5-i)} 
$\boldsymbol{\Lambda} (P) =\langle \nots,{\bf A},{\bf D},{\bf B},\noti\rangle$.
\newline 
{\rm (5-ii)}  $PR_P$ is  principal, $P$ verifies $(\boldsymbol{\natural})$ and is the only divisorial $P$-primary ideal. \newline 
{\rm (5-iii)}   $PR_P$ is  principal,  each finitely generated ideal contained in $P$ is contained in at least one maximal ideal that does not contain $P$ and $P$  verifies $(\boldsymbol{\natural})$. \newline 
{\rm (5-iv)}   $PR_P$ is  principal,   each finitely generated ideal contained in $P$ is contained in infinitely many maximal ideals that do not contain $P$ and $P$ verifies $(\boldsymbol{\natural})$.

\item[(6)] The following are equivalent: \\ 
{\rm (6-i)} $\boldsymbol{\Lambda} (P) = \langle \nots,{\bf A},{\bf D},\notb, {\bf I}\rangle$. \newline
{\rm (6-ii)}  $P$ is  the only $P$-primary ideal, each finitely generated ideal contained in $P$ is contained in at least one maximal ideal that does not contain $P$ and $P$ verifies $(\boldsymbol{\natural})$. \newline 
{\rm (6-iii)}   $P$ is  the only $P$-primary ideal, each finitely generated ideal contained in $P$ is contained in infinitely many maximal ideals that do not contain $P$ and $P$ verifies $(\boldsymbol{\natural})$. 

 \item[(7)] The following are equivalent: \\
 {\rm (7-i)} 
$\boldsymbol{\Lambda} (P) =\langle \nots, \nota, {\bf D},{\bf B},{\bf I}\rangle$. \newline
{\rm (7-ii)}  $PR_P$ is not  principal, proper $P$-primary ideals exist, but $P$ is the only divisorial $P$-primary ideal and $P$ does not verify $(\boldsymbol{\natural})$. \newline 
{\rm (7-iii)}  $PR_P$ is not  principal, $P$ is divisorial and minimal over some finitely generated ideal, but each finitely generated ideal contained in $P$ is contained in at least one maximal ideal that does not contain $P$ and $P$ does not verify $(\boldsymbol{\natural})$. \newline 
{\rm (7-iv)}  $PR_P$ is not principal, $P$ is divisorial and minimal over some finitely generated ideal, but each finitely generated ideal contained in $P$ is contained infinitely many maximal ideals that do not contain $P$ and $P$ does not verify $(\boldsymbol{\natural})$.

\item[(8)] The following are equivalent: \\
{\rm (8-i)} 
$\boldsymbol{\Lambda} (P) =\langle \nots,\nota,{\bf D},{\bf B},\noti\rangle$.
\newline 
{\rm (8-ii)}  $PR_P$ is  principal, $P$ is the only divisorial $P$-primary ideal and $P$ does not verify $(\boldsymbol{\natural})$. \newline 
{\rm (8-iii)}   $PR_P$ is  principal, $P$ is divisorial, each finitely generated ideal contained in $P$ is contained in at least one maximal ideal that does not contain $P$and $P$ does not verify $(\boldsymbol{\natural})$. \newline 
{\rm (8-iv)}   $PR_P$ is  principal, $P$ is divisorial, each finitely generated ideal contained in $P$ is contained in infinitely many maximal ideals that do not contain $P$ and $P$ does not verify $(\boldsymbol{\natural})$.

\item[(9)] The following are equivalent: \\ 
{\rm (9-i)}  $\boldsymbol{\Lambda} (P) = \langle \nots, \nota,{\bf D},\notb,{\bf I}\rangle$. \newline
{\rm (9-ii)}  
$P$ is divisorial and the only $P$-primary ideal, each finitely generated ideal contained in $P$ is contained in at least one maximal ideal that does not contain $P$ and $P$ does not verify $(\boldsymbol{\natural})$. \newline 
{\rm (9-iii)}   $P$ is divisorial and the only $P$-primary ideal, each finitely generated ideal contained in $P$ is contained in infinitely many maximal ideals that do not contain $P$ and $P$ does not verify $(\boldsymbol{\natural})$.

\item[(10)]  The following are equivalent:  \\
{\rm (10-i)} 
$\boldsymbol{\Lambda} (P) = \langle \nots, \nota,\notd,{\bf B},{\bf I}\rangle$. 
\newline
{\rm (10-ii)}   $PR_P$ is not principal, $P$ is not divisorial and it is minimal over some finitely generated ideal, but each finitely generated ideal contained in $P$ is contained in at least one maximal ideal that does not contain $P$. \newline
{\rm (10-iii)}  $PR_P$ is not principal, $P$ is not divisorial and it is minimal over some finitely generated ideal, but each finitely generated ideal contained in $P$ is contained infinitely many maximal ideals that do not contain $P$.

\item[(11)]  The following are equivalent:  \\
{\rm (11-i)} 
$\boldsymbol{\Lambda} (P) = \langle \nots, \nota,\notd,{\bf B},\noti,\rangle$.
\newline
{\rm (11-ii)}   $PR_P$ is principal,  $P$ is not divisorial and each finitely generated ideal $I\subseteq P$ is contained in at least one maximal ideal that does not contain $P$. \\
{\rm (11-iii)}   $PR_P$ is principal,  $P$ is not divisorial and each finitely generated ideal $I\subseteq P$ is contained in infinitely many maximal ideals that do not contain $P$. \\
{\rm (11-iv)}   $PR_P$ is principal,  $P$ is not divisorial and for each finitely generated ideal $I\subseteq P$, there is an ideal containing $I$ that is incomparable with $P$.

\item[(12)]  The following are equivalent:  \\
{\rm (12-i)}   $\boldsymbol{\Lambda} (P) = \langle \nots, \nota,\notd,\notb,{\bf I} \rangle$. \\
 {\rm (12-ii)}  $P$ is the only $P$-primary ideal, $P$ is not divisorial and each  finitely generated ideal contained in $P$ is contained in at least one maximal ideal that does not contain $P$. \\
{\rm (12-iii)}  $P$ is the only $P$-primary ideal, $P$ is not divisorial and each  finitely generated ideal contained in $P$ is contained in infinitely many  maximal ideals that do not contain $P$.

\end{enumerate} \end{theorem}

\begin{proof}   Let   $P$ be nonzero, nonmaximal ideal of $R$.  If  $P$ is  sharp then (by Remark \ref{rk:1.1}, or Corollary \ref{cor:1.4})  it is necessarily divisorial. Thus the eight cases $\langle {\bf S}, \mbox{--}, \notd, \mbox{--}, \mbox{--} \rangle$,  arising from any  values given to $\mbox{--}$   are impossible.  Moreover, if $P$ is sharp, then  it is also antesharp by Corollary \ref{cor:1.4} so  the  four cases $\langle {\bf S},\nota, {\bf D}, \mbox{--}, \mbox{--} \rangle$,  arising from any  values given to  $\mbox{--}$   are impossible.  As with maximal ideals, clearly  $\notb$ implies ${\bf I}$  \cite[Theorem 23.3(e)]{Gilmit}.  So $\boldsymbol{\Lambda}(P)$ cannot contain both $\noti$ and $\notb$, eliminating  five more cases of the form $\langle  \mbox{--}, \mbox{--},  \mbox{--}, \notb,\noti, \rangle$.   Finally, we observed in Corollary \ref{cor:1.4} that a nonzero nonmaximal  antesharp prime ideal  is necessarily divisorial cutting three more cases. The remaining twelve (out of thirty-two) cases are those listed in the statement.
 
As in the maximal case, $P$ is idempotent if and only if  $PR_P$ is not  principal. Also as mentioned above, if it is unbranched, it must be idempotent. Thus $PR_P$ principal is equivalent to $P$ is  branched and not idempotent. Also sharp implies antesharp and antesharp implies divisorial by Corollary \ref{cor:1.4}. 

  (1) By Proposition \ref{pr:1.7} [respectively, Proposition \ref{pr:1.9}] (1-ii) and (1-iii) [respectively, (1-ii) and (1-iv)] are equivalent. It is clear that (1-i) implies (1-iii). Finally, (1-ii) implies (1-i)
by Remark \ref{rk:1.1}, Proposition \ref{pr:1.9} and by the fact that  sharp implies antesharp.

(2) As noted above, the condition ``$PR_P$ is principal" is equivalent to  ``$P$ is branched and not idempotent".    Corollary \ref{cor:1.4} and Proposition \ref{pr:1.9} are the only other things needed  for the equivalence of (2-i) through (2-iv). Since $R$ is Pr\"ufer, $P=P(P:P)$ is a prime ideal of $(P:P)$ and all other primes of $(P:P)$ are of the form $Q'=Q(P:P)$ for some prime $Q$ of $R$ with $R_Q=(P:P)_{Q'}$. In particular $R_P=(P:P)_P$ and if $N$ is a maximal ideal of $R$ that does not contain $P$, then $N'=N(P:P)$ is a maximal ideal of $(P:P)$ that does not contain $P$ and  $R_N=(P:P)_{N'}$. If $P$ is invertible as an ideal of $(P:P)$, then it is maximal, divisorial, sharp, branched and locally principal as a prime of $(P:P)$. Hence $R_P=(P:P)_P$ does not contain $\Theta_R(P)=\Theta_{(P:P)}(P)$.  Thus  $P$ is sharp, antesharp, divisorial, branched and not idempotent as a prime of $R$. Conversely, if $P$ is locally principal and the radical of a finitely generated ideal, then there is a finitely generated ideal $B\subsetneq P$ with $\sqrt B=P$ and $BR_P=PR_P$. Since we have already established that $P$ is antesharp in  this situation, checking locally in $(P:P)$ shows that $P=B(P:P)$ is an invertible ideal of $(P:P)$.

(3) The equivalence of these three statements is by Proposition \ref{pr:1.2}, Corollary \ref{cor:1.4}   and the fact that an unbranched prime is idempotent. 

Suppose $I$ is a finitely generated ideal contained in $P$ such that $M_1$, $M_2$, \dots , $M_n$ are the only maximal ideals that contain $I$ and do not contain $P$. By choosing elements $a_i\in P\backslash M_i$, we have a finitely generated ideal $J:=I+(a_1,\dots,a_n)\subsetneq P$ such that each maximal ideal that contains $J$ also contains $P$. Hence $P$ is sharp. So having $P$ not sharp is also  equivalent to each finitely generated ideal contained in $P$ is contained in infinitely many maximal ideals that do contain $P$. We use this equivalence in all of the remaining cases.

(4) For this case, we again use that $P$ is both idempotent and branched if and only if $PR_P$ is not principal and $P$ is minimal over a finitely generated ideal. Proposition \ref{pr:1.3} and Corollary \ref{cor:1.4} take care of establishing the equivalence of  ``$P$ is divisorial and antesharp" and   ``$P$ verifies $(\boldsymbol{\natural})$". Use either Proposition \ref{pr:1.7} or Proposition \ref{pr:1.9} and the characterization of ``not sharp" above to finish this case. 

 (5) As in (2), $PR_P$ is principal if and only if $P$ is branched and not idempotent. Finish by using the rest of the argument for case (4).

(6) The definition of an unbranched prime  together with the last two sentences from the proof of  (4) are all that is needed here.

(7, 8, and 9) By  Proposition \ref{pr:1.3}, $P$ is not antesharp  if and only if it does not verify $(\boldsymbol{\natural})$. Thus Corollary \ref{cor:1.4} no longer applies, and we need to include the assumption $P$ is divisorial in (7-iii), (7-iv), (8-ii), (8-iii), (8-iv), (9-ii) and (9-iii). Other than that, the proof for (7) is essentially the same as for (4), the one for (8) matches up with (5) and the one for (9) with (6).

(10, 11 and 12) Match (10) with (7) and (4), (11) with (8) and (5), and (12) with (9) and (6). Then replace ``$P$ is divisorial" by ``$P$ not divisorial" and finish using appropriate parts of the respective proofs for (4), (5) and (6). 
 \end{proof}

\section{Examples}\label{Examples}

The goal of the present section is to show, by giving explicit examples, that all the cases considered in Theorem \ref{th:1.11} and Theorem \ref{th:1.12} are attained.

Let $R$ be an integral domain  and $X$ an indeterminate over $R$. For each polynomial $f \in R[X]$, let $\boldsymbol{c}(f)$ denote the content of $f$ (i.e. the ideal of $R$ generated by the coefficients of $f$). Recall that the Nagata ring $R(X)$ is defined as follows:
$$
R(X) := \{  f/g   \mid f, g \in R[X] \mbox{ and }  \boldsymbol{c}(g)= R \}  \mbox{ (cf. \cite[page 18]{Nagata:1962} and \cite[Section 33]{Gilmit})}.
$$

\begin{theorem} \label{th:2.1} Let $P$ be a nonzero prime of a Pr\"ufer domain $R$. Then $PR(X)$ is a prime in the Pr\"ufer domain $R(X)$ and $P$ and $PR(X)$ have the same type, i.e.:
$$
\boldsymbol{\Lambda}_R(P) =  \boldsymbol{\Lambda}_{R(X)}(P(X)).$$
 \end{theorem}
 
 We start by proving the following lemma which is likely known, but we are unable to find a reference (T. Nishimura  proved a similar statement for the $v$-operation for the ring of polynomials $R[X]$ with coefficients in an integral domain $R$ \cite[Proposition 7.1]{Nishimura:1961} and the third named author  extended this result to the case of rings with zero divisors \cite[Lemma 5.1]{Lukrull}).

 \begin{lemma}\label{le:2.2} Let $R$ be an integral domain. If $I$ is a nonzero ideal of $R$, then $I$ is divisorial  in $R$ if and only if $IR(X)$ is divisorial in $R(X)$.\end{lemma}

\begin{proof} We first show that $(R(X) : IR(X))=(R : I)R(X)$ (even if $I$ is not finitely generated).   Clearly, if $t\in (R:I)$, then $tIR(X)\subseteq R(X)$. So $(R(X):IR(X))\supseteq (R:I)R(X)$. For the reverse containment, let $z\in (R(X):IR(X))$ and let $i$ be a nonzero element of $I$. Then there is an element  $s \in R(X)$ such that $z=s/i$. Also there is a polynomial $g\in R[X]$ with unit content such that $f:=gs\in R[X]$. So
$gz=f/i\in K[X]$. For each $j \in I$, $gzj=fj/i \in R(X) \cap K[X] =R[X]$. It follows that,  for each coefficient $b_k$ of $f$, $(b_k/i)j \in R$. So $f/i\in (R:I)R[X]$ and we have $ z =f/ig\in (R:I)R(X)$ as desired. 

From the previous property, it follows immediately that, if  $H$ is a nonzero fractional ideal of $R$, then  $(R(X) : HR(X))=(R : H)R(X)$. The result follows. 
  \end{proof}
 
 \smallskip 

\begin{proof} (\bf Theorem \ref{th:2.1}\rm ) By Lemma \ref{le:2.2}, $P$ is a divisorial ideal of $R$  if and only if $PR(X)$ is a divisorial ideal of $R(X)$.

 It is easy to see that $P^2R(X)=(PR(X))^2$, thus $P$ is idempotent if and only if $PR(X)$ is idempotent.  
 
 It is well-known that if $Q$ is a $P$--primary ideal of $R$ then $Q(X)$ is a $P(X)$--primary ideal of $R(X)$  \cite[(6.17)]{Nagata:1962} and, conversely, it is straightforward to verify that if $H$ is a  $P(X)$--primary ideal of $R(X)$   then $Q:=H \cap R$ is a $P$--primary ideal of $R$. Therefore $P$ is branched if and only if $PR(X)$ is branched. 
 
It is clear that if $P$ is sharp, then $PR(X)$ is sharp.  Conversely, let  $PR(X)$ be a sharp prime ideal of $R(X)$ and let  $H =(z_1,z_2,\dots, z_n)\subseteq P(X)$ be a finitely generated ideal of $R(X)$ such that the only maximal ideals of $R(X)$  that contain $H$ are those that contain $P(X)$. Then, for each $k$, $1 \leq k \leq n$, there exists a polynomial with unit content $g\in R[X]$ such that $f_k := z_kg\in P[X]$. Let $a_{k,j}$, with $0\leq j \leq t_k$, be the coefficients of the polynomial $f_k$ (of degree $t_k$). Set $I := (a_{k, j} \mid 1\leq k \leq n, 0\leq j \leq t_k)$. Then clearly $I$ is a finitely generated ideal of $R$ contained in $P$. If $I$ was contained in some maximal ideal $N$, with $P \not\subseteq N$, then clearly $f_k \in N[X]$ and so $z_k \in N(X)$, for each $k$, with $N(X)$ maximal ideal of $R(X)$,  with $P(X) \not\subseteq N(X)$ and this is impossible since $H \not\subseteq N(X)$ by assumption. Therefore if $P(X)$ is sharp, then $P$ is also sharp.

Note that if $P$ is a nonzero nonmaximal  prime ideal of $R$ then $P^{-1} = (P:P)$ and $P(X)$ is a nonzero nonmaximal  prime ideal of $R(X)$,  hence (by the proof of Lemma \ref{le:2.2}) we have $((P(X):P(X))= P(X)^{-1} = (R(X):P(X)) = (R:P)(X) =P^{-1}(X) = (P:P)(X)$. Thus $P$ is maximal in $(P:P)$ (i.e. $P$ is antesharp) if and only if $P(X)$ is maximal in $((P(X):P(X))$, (i.e. $P(X)$ is antesharp).\end{proof}

\begin{theorem} \label{th:2.3} Let $M$ be a nonzero maximal ideal of a Pr\"ufer domain $T$ with corresponding residue field $K$. Let $\varphi: T(X) \rightarrow T(X)/M(X) \cong K(X)  $ be the canonical projection and let $R$ be the pullback of $K[X]_{(X)}$ along $M(X)$ (i.e. $R:= \varphi^{-1} (K[X]_{(X)})$), then $R$ is a Pr\"ufer domain such that $M(X)$ is a nonmaximal divisorial prime of $R$. Also $M(X)$ retains its other four properties, i.e.
$\boldsymbol{\Lambda}_R (M(X)) = \langle{\bf V},{\bf A},{\bf D},{\bf Y,Z}\rangle$ where $\boldsymbol{\Lambda}_{T(X)} (M(X)) =\langle{\bf V,X,Y,Z}\rangle (= \boldsymbol{\Lambda}_{T} (M)$ by Theorem \ref{th:2.1}).
 \end{theorem}

\begin{proof} Since $M(X)$ is a maximal ideal of $T(X)$, it is the conductor of $R$ into $T(X)$. Thus it is a divisorial prime of $R$ that is not maximal (note that in $R$, $M(X) = \varphi^{-1}(0) \subsetneq \varphi^{-1}(XK[X]_{(X)})=: \mathfrak{M})$. By  pullback properties \cite{Fontop}, each maximal ideal of $R$ that does not contain $M(X)$ is the contraction of a unique maximal ideal of $T(X)$ and $\mathfrak{M}$ is the only maximal ideal of $R$ that contains $M(X)$.  Thus it is easy to see that $M(X)$ is sharp [respectively:  idempotent, branched, antesharp] in $R$ if and only if $M(X)$ is sharp  [respectively:  idempotent, branched, antesharp]  in $T(X)$.  Note that since $M$ is a maximal ideal of $T$, $M$ is  obviously antesharp in $T$  and so $M(X)$ is antesharp in $R$. \end{proof}

Since each nonzero nonmaximal prime $P$ of a valuation domain is sharp and divisorial  and the (nonzero) maximal ideal  $M$ of a valuation domain is always sharp, the valuation domains give 
the first easy  examples for  special values assumed by $\boldsymbol{\Lambda}(M)$ and $\boldsymbol{\Lambda}(P)$. In all of the examples, $M$ will be used to denote the specific maximal ideal of the desired type and $P$ will be used to denote the specific prime ideal of the desired type. For some pairs of examples,  the specific prime $P$ represents in the earlier example may appear in the later one with a new name.  The pullback construction in Theorem \ref{th:2.3} is used in several examples, specifically in Examples 4.8--4.11. A related construction is used in Examples 4.17 and 4.19.

\begin{example} $\boldsymbol{\Lambda}(M) = \langle {\bf S},{\bf D}, {\bf B},\noti\rangle$  (Case 1  of Theorem \ref{th:1.11}).

Take the maximal ideal $M$ of a discrete valuation domain.
 \end{example}

\begin{example}
$\boldsymbol{\Lambda}(M)=\langle {\bf S},\notd, {\bf B}, {\bf I}\rangle$  (Case 2  of Theorem \ref{th:1.11}).

  Take  $V$ a one-dimensional valuation domain, not a DVR. In this case, $M$ is  both branched and idempotent but it is not  divisorial. 
 \end{example}

\begin{example} $\boldsymbol{\Lambda}(P) = \langle {\bf S},{\bf A},{\bf D},{\bf B}, {\bf I} \rangle$ (Case 1  of Theorem \ref{th:1.12}).  

 Take a two-dimensional valuation domain $V$ with height-one prime $P$ such that $V_P$ is not a DVR. 
 \end{example}

\begin{example} $\boldsymbol{\Lambda}(P) = \langle {\bf S},{\bf A},{\bf D},{\bf B},\noti\rangle$ (Case 2  of Theorem \ref{th:1.12}).   

Take a two-dimensional valuation domain $V$ with height-one prime $P$ such that $V_P$ is  a DVR. \end{example}

\begin{example}  $\boldsymbol{\Lambda}(M) = \langle {\bf S},\notd,\notb,{\bf I} \rangle$ (Case 3  of Theorem \ref{th:1.11})  \ and \\ $\mbox{ }$  \hskip 2.6cm$\boldsymbol{\Lambda}(P) = \langle {\bf S},{\bf A},{\bf D},\notb,{\bf I}\rangle$ (Case 3  of Theorem \ref{th:1.12}).  

Start with a valuation domain $V$ where the maximal ideal $M$  is unbranched (hence, in particular, $M$ is not divisorial) then the same is true in $V(X)$. We have $\boldsymbol{\Lambda}(M)=\langle {\bf S}, \notd,\notb,{\bf I}\rangle$. Let $K$ be the residue field $V/M$, let $\varphi: V(X) \rightarrow V(X)/M(X) \cong K(X)  $ be the canonical projection and let $R$ be the pullback of $K[X]_{(X)}$ along $M(X)$ (i.e. $R:= \varphi^{-1} (K[X]_{(X)})$). It is easy to see that $R$ is a valuation domain with $P:=M(X)$ a divisorial nonmaximal prime that is unbranched and sharp (thus, antesharp) (Theorem \ref{th:2.3}). Thus  $\boldsymbol{\Lambda}_R(P) = \langle {\bf S},{\bf A},{\bf D},\notb,{\bf I}\rangle$. \end{example}

\begin{example} $\boldsymbol{\Lambda}(M) = 
\langle \nots,\notd,{\bf B},{\bf I} \rangle$  (Case 4  of Theorem \ref{th:1.11})   \ and \\ $\mbox{ }$  \hskip 2.6cm
$\boldsymbol{\Lambda}(P) = 
\langle \nots,{\bf A},{\bf D},{\bf B},{\bf I}\rangle$  (Case 4  of Theorem \ref{th:1.12}). 

Let $T$ be the ring of algebraic integers and let $M$ be a maximal ideal of $T$. Then $M$ is branched and idempotent ($M$ has these properties, since $T_M$  is a one-dimensional nondiscrete valuation domain), but $M$ is not sharp since each nonunit is contained in infinitely many maximal ideals.  Thus  $\boldsymbol{\Lambda}(M) = 
\langle \nots,\notd,{\bf B},{\bf I} \rangle$.  Let $K$ be the residue field $T/M$ and let $\varphi: T(X) \rightarrow T(X)/M(X) \cong K(X)  $ be the canonical projection and let $R$ be the pullback of $K[X]_{(X)}$ along $M(X)$. Then, in $R$, the nonmaximal prime  ideal $P:=M(X)$ is not sharp but it is antesharp, divisorial, branched and idempotent  (Theorem \ref{th:2.3}). Hence  $\boldsymbol{\Lambda}_R(P)=\langle \nots,{\bf A},{\bf D},{\bf B},{\bf I}\rangle$.
 \end{example}

  \begin{example}
$\boldsymbol{\Lambda}(M) = \langle \nots,\notd,{\bf B},\noti \rangle$  (Case 5  of Theorem \ref{th:1.11})   \ and \\ $\mbox{ }$  \hskip 2.85cm $\boldsymbol{\Lambda}(P) = 
\langle \nots,{\bf A},{\bf D},{\bf B},\noti\rangle$ (Case 5  of Theorem \ref{th:1.12}).

 Let $D$ be an almost Dedekind domain that is not Dedekind (for instance \cite[Example 8.2]{FHPbook}) and let $M$ be a maximal ideal of $D$ that is not invertible (or, equivalently, not divisorial). Thus 
$\boldsymbol{\Lambda}(M) = \langle \nots,\notd,{\bf B},\noti \rangle$.  The ring  $D(X)$ is also an almost Dedekind domain \cite[Proposition 36.7]{Gilmit}  with $M(X)$ maximal and not divisorial. So $M(X)$ is branched but not sharp  in $D(X)$. As above,  let $K$ be the residue field $D/M$ and let $\varphi: D(X) \rightarrow D(X)/M(X) \cong K(X)  $ be the canonical projection and let $R$ be the pullback of $K[X]_{(X)}$ along $M(X)$. Then $P:= M(X)$ is a nonmaximal prime ideal  in $R$, which is neither sharp nor idempotent, but it is antesharp, divisorial and branched   (Theorem \ref{th:2.3}). Thus $\boldsymbol{\Lambda}_R(P)=\langle \nots,{\bf A},{\bf D},{\bf B},\noti\rangle$. \end{example}

\begin{example}  $\boldsymbol{\Lambda}(M) = \langle \nots,\notd,\notb,{\bf I} \rangle$ (Case 6  of Theorem \ref{th:1.11})  \ and \\ $\mbox{ }$  \hskip 2.85cm 
  $\boldsymbol{\Lambda}(P) = 
\langle \nots,{\bf A},{\bf D},\notb,{\bf I}\rangle$ (Case 6  of Theorem \ref{th:1.12}).

 Let $T$ be the Pr\"ufer domain constructed in \cite[Example 38]{Lurtp} and let $M$ be the unbranched maximal ideal of $T$. In $T$, $M$ is  not divisorial and not sharp since each finitely generated
ideal contained in $M$ is contained in infinitely many other maximal ideals of $T$. Thus  $\boldsymbol{\Lambda}(M) = \langle \nots,\notd,\notb,{\bf I} \rangle$. As above,  let $K$ be the residue field $T/M$ and let $\varphi: T(X) \rightarrow T(X)/M(X) \cong K(X)  $ be the canonical projection and let $R$ be the pullback of $K[X]_{(X)}$ along $M(X)$.  Then $R$ is a Pr\"ufer domain where $P:= M(X)$ is a divisorial nonmaximal prime. In $R$ it is  antesharp and unbranched (and idempotent) but not sharp   (Theorem \ref{th:2.3}) so $\boldsymbol{\Lambda}_R(P)=\langle \nots,{\bf A},{\bf D},\notb,{\bf I}\rangle$.  
\end{example} 

 Let $D$ be a Pr\"ufer domain with a nonmaximal prime $Q$ such that $Q^{-1}=D$ and let $X$ be an indeterminate over $D$. Then, in the Nagata ring  $S:=D(X)$, the ideal $P:=Q(X)$ is a nonmaximal prime with $P^{-1} \ (=Q(X)^{-1})=S$. Let $K:=D_Q/QD_Q$. Then $S_P/PS_P$ is naturally isomorphic to $K(X)$.  Let $\psi_P: S_P \rightarrow S_P/PS_P \cong K(X)$ be the canonical projection,  let $W:=K[X]_{(X)}$  and let $V$ be the pullback of $W$ along $PS_P$ (i.e. $V:= \psi_P^{-1} (W)$).

\begin{theorem} \label{th:2.12} Let $D$, $Q$, $S$,  $P$, $W$, $V$ be as above and let  $R:=V\cap S$. 
\begin{enumerate}

\item If $ \psi: S \rightarrow K(X) $ is the canonical homomorphism (obtained by composition of the natural inclusion $j: S \rightarrow S_P$ with the canonical projection  $\psi_P: S_P \rightarrow S_P/PS_P \cong K(X)$), then $R = \psi^{-1} (W)$.

\item $P$ is a prime ideal of $R$ and $(P:P)=(R:P)=S$.

\item $XR$ is  a principal maximal ideal of $R$ that  properly contains $P$. Moreover, $\bigcap_{n\geq 1} X^nR=P$. 

\item $R$ is a B\'ezout domain.

\item $P$ is a divisorial nonmaximal prime ideal of $R$ that is not antesharp. 
\end{enumerate}
\end{theorem}

\begin{proof}  (1)  By  well-known properties of pullbacks, $V$ is a valuation domain having $S_P$ as an overring and $PS_P$ as a prime     nonmaximal ideal and clearly $R = S \cap V = j^{-1}(V)  = 
j^{-1}(\psi_P^{-1} (W)) = \psi^{-1} (W)$. 

(2)  Since $PS_P\subsetneq V$, $P \ (= PS_P  \cap S = PS_P \cap R)$ is a prime (but not maximal) ideal of $R$. By the proof of Lemma \ref{le:2.2}, $(S:P)=(D(X):Q(X))=(D:Q)D(X)=D(X)=S$. Thus  $S \subseteq (P:P)\subseteq (R:P)  \subseteq (S:P)  = S $.

(3) 
Set $N:=XV \cap S$. Then  $(1/X)N\subseteq V$ and $(1/X)N\subseteq S$ since $X$ is a unit of $S$. It follows that $N=XR$. Also, for each nonunit $a\in R\backslash N$,  the element $a+X^k$ is a unit of $V$ for each positive integer $k$ (since $V$ is a valuation domain with maximal ideal $XV$  and $a \in V \setminus XV$). On the other hand $a \in S=D(X)$, thus we may write $a=b/u$ where $b,u\in D[X]$ with $\boldsymbol{c}(u)=D$. Clearly, for sufficiently large $m$, $b+uX^m$ also has unit content in $D$ and therefore, $(b+uX^m)/u =a+X^m$ is a unit in $D(X)=S$ and so $a+X^m$ is a unit in $R$, for sufficiently large $m$. It follows that $N$ is a maximal ideal of $R$.  

Moreover, note that $\bigcap_{n\geq 1} X^nW=0$ and so, by the properties of the pullback constructions, the maximal ideal $XV \ (=\psi_P^{-1} (XW))$  of the valuation domain $V \ (=\psi_P^{-1} (W))$ is such that  $\bigcap_{n\geq 1} X^nV = \psi_P^{-1} (0) =PS_P$ and so $\bigcap_{n\geq 1} X^nR=P$.

(4) Let $A=(a_0,a_1,\dots,a_n)$ be a finitely generated (nonzero) ideal of $R$. Without loss of generality, we may assume $AV=a_0V$. Each $a_i$ has the form $b_i/u_i$ with $b_i,u_i\in D[X]$ and $\boldsymbol{c}(u_i)=D$. There is no harm in assuming $u=u_i$ for each $i$.  Recursively define integers $k_0: =0$ and  $k_i:=\mbox{deg(}b_{i-1})+k_{i -1}+1$. Then set $b:=\sum_ i b_iX^{k_i} \in D[X]$. It is easy to see that, in $D$,   $ \boldsymbol{ c}(b) =\sum_i \boldsymbol{ c}(b_i)$.  Consider the element $a:=b/u \ (\in D(X))$. Clearly, $a \in A$ and from our assumption that $AV=a_0V$, we have that $a_iX^{k_i}V$ is properly contained in $AV$ for each positive integer $i$. Thus $AV=aV$.  We also have $aS= \boldsymbol{c}(b)S=AS$. Thus, as with $N$, we must have $A=aR$. Hence $R$ is a B\'ezout domain. 

(5) That $P$ is a divisorial (prime) ideal of $R$ follows from the fact that $P=\bigcap_{n\geq 1} X^nR$. On the other hand, $P$ is not antesharp in $R$ since $(P:P)=S$ where $P$ is not maximal.    \end{proof}
 
The assumption in Theorem \ref{th:2.12} that $Q$ has a trivial dual is  a little stronger than one needs to get a divisorial prime that is not antesharp. In fact, one can start with an arbitrary nonempty set of nonzero nondivisorial primes and build a Pr\"ufer domain with a corresponding set of divisorial primes that are not antesharp whenever the corresponding contraction to the original ring is not maximal.  To set the notation for the next theorem let $\mathcal Q:=\{Q_\alpha \mid \alpha \in \mathcal A\}$ be a set of nonzero nondivisorial prime ideals of a Pr\"ufer domain $D$ with  $|\mathcal A|\ge 2$. Next let $\mathcal Y:=\{Y_\alpha \mid \alpha \in \mathcal A\}$ be a set of algebraically independent indeterminates over $D$ and set $S:=D(\mathcal Y)$. For each $\alpha\in \mathcal A$, let $P_\alpha =Q_\alpha(\mathcal Y)$, $K_\alpha =D_{Q_\alpha}/Q_\alpha D_{Q_\alpha}$ and $\mathcal Y_\alpha=\mathcal Y\backslash \{Y_\alpha\}$. Then the quotient field of $S/P_\alpha$ is naturally isomorphic to $K_\alpha(\mathcal  Y)$. Let $V_\alpha$ be the pullback of $W_\alpha :=K_\alpha(\mathcal Y_{\alpha})[Y_\alpha]_{( Y_\alpha)}$ along  $P_\alpha S_{P_\alpha}$ (i.e.  if $\psi_{P_\alpha}:  S_{P_\alpha} \rightarrow S_{P_\alpha}/P_{\alpha} S_{P_\alpha} \cong K_\alpha(\mathcal  Y)$ is the canonical projection,  then $V_{\alpha} := (\psi_{P_\alpha})^{-1} (W_\alpha)$).  Finally, let $R: = (\bigcap_{\alpha \in \mathcal A}  V_\alpha) \cap S$ and let $P'_\alpha:=
P_\alpha\cap R$ for each $\alpha\in \mathcal A$.

\begin{theorem} \label{th:2.13} Let $D$, $\mathcal Q$, $\mathcal Y$, $\{V_\alpha\}$, $R$, $S$ and $\{P'_\alpha\}$ be as above and let    $\mathcal M:=\{P_\alpha\mid Q_\alpha\in {\rm Max}(D)\}$.  
\begin{enumerate}
\item $R$ has the same quotient field as $S$.
\item $N_\alpha:=Y_\alpha V_\alpha \cap R$ is a principal maximal ideal of $R$ generated by $Y_\alpha$. 
\item Each nonzero two-generated ideal of $R$ is principal.
\item $R$ is a B\'ezout domain.
\item ${\rm Max}(R)=\{M \cap R \mid M\in \mbox{\rm Max}(S)\backslash \mathcal M\}  \cup \{N_\alpha \mid \alpha \in \mathcal A\}$.  
\item $(R:P'_\alpha)= (\bigcap \{V_\gamma \mid \gamma \in \mathcal A \mbox{ and }  N_\gamma+ P'_\alpha=R\}) \cap (S:P_\alpha)$. 
\item Each  $P'_\alpha$ is a divisorial prime ideal of $R$, and $P'_\alpha$ is  antesharp as  ideal of $R$ if and only if $P_\alpha\in \mathcal M$  (equivalently, $Q_\alpha$ is a maximal ideal of $D$). \end{enumerate} \end{theorem}

\begin{proof} {(1)}  Note that each $V_\alpha$ contains $D[\mathcal Y]$. Thus $R$ and $S$ have the same quotient field. For each $V_\alpha$, we let $v_\alpha$ be the corresponding valuation map. 

(2) It is clear that  $(1/Y_\alpha)N_\alpha\subseteq V_\alpha$. Since $Y_\alpha$ is a unit in all other $V_\beta$'s and in $S$, we also have $(1/Y_\alpha)N_\alpha \subseteq V_\beta$ for each $\beta$ and $(1/Y_\alpha)N_\alpha \subseteq S$. It follows that $N_\alpha =Y_\alpha R$. Clearly, $N_\alpha$ is a prime ideal of $R$. 

(3 and 4)  Let $A$ be a two-generated ideal of $R$. Then without loss of generality we may assume there are polynomials $a,b, u\in D[\mathcal Y]$ with ${\bf c}(u)=D$ such that $A=(a/u, b/u)$. Let $m$ be a positive integer that is larger than the sum of the total degree of $b$ and the total degree of $a$. Next select two distinct members of $\mathcal Y$, say $Y_\beta$ and $Y_\gamma$. Set $t=Y_\beta^m+Y_\gamma^{2m+1}$ and $d=at+b$.  Since there is no overlap in degrees for the individual terms in $at$ and $b$ and at least one of $Y_\beta$ and $Y_\gamma$ is a unit in $V_\alpha$, the value of $v_\alpha(d/u)$ is the minimum of $v_\alpha(a/u)$ and  $v_\alpha(b/u)$. Thus $(d/u)V_\alpha=AV_\alpha$. Since $u$ is a unit of $S$ and there is no overlap in degrees for the individual terms in $at$ and $b$, $(d/u)S=dS\subseteq AS=(a,b)S\subseteq {\bf c}(d)S$. Since  $D$ is Pr\"ufer, there is a polynomial $f\in {\bf c}(d)^{-1}[\mathcal Y]$ such that $df$ has unit content in $D$ and $f{\bf c}(d)S\subseteq S$. Hence $df$ is  a unit of $S$ and we have   ${\bf c}(d)S=df{\bf c}(d)S\subseteq (d/u)S\subseteq AS\subseteq {\bf c}(d)S$. Thus as with $N_\alpha$, $A$ is principal. Therefore, $R$ is B\'ezout and each $N_\alpha$ is a principal maximal ideal of $R$. 

(5) Next we  show  ${\rm Max}(R)= \{M \cap R \mid M\in \mbox{\rm Max}(S)\backslash \mathcal M\}  \cup \{N_\alpha \mid \alpha \in \mathcal A\}$. Clearly, if $Q_\alpha$ is  a maximal ideal of $D$, then $P_\alpha$ is a maximal ideal of $S$, but $P'_\alpha$ is not a maximal ideal of $R$ -- it is properly contained in $N_\alpha$. Let $N$ be a maximal ideal of $R$. There is nothing to prove if $N=N_\alpha$ for some $\alpha$.  To see that $N$ is the contraction of a maximal ideal of $S$, it suffices to show that $R_N$ contains $S$. Let $u\in D[\mathcal Y]$ be a polynomial with unit content in $D$. Then $u$ is a unit of $S$ and it is unit in $V_\alpha$ if it is not a multiple of $Y_\alpha$ in $D[\mathcal Y]$. Since $u$ is a polynomial there are at most finitely many $\alpha_i$'s such that $Y_{\alpha_i}$ divides $u$ in $D[\mathcal Y]$. Hence we may write $u=Y^{n_1}_{\alpha_1}\cdots Y^{n_m}_{\alpha_m}w$ where $w$ has unit content in $D$ and no $Y_\alpha$ divides $w$ in $D[\mathcal Y]$. Thus $w$ is unit in $S$ and in each $V_\alpha$ which makes it a unit of $R$. Since each $N_\alpha=Y_\alpha R$ is a principal maximal,  no $Y_\alpha$ is contained in $N$. It follows that $u$ is a unit of $R_N$ and therefore, $R_N$ contains $S$.

(6) Now consider what happens to $P'_\alpha=P_\alpha\cap R$ as a prime of $R$. Since $Q_\alpha\subsetneq  D\subsetneq  V_\beta$ for each $\alpha$ and $\beta$, $P'_\alpha$ contains $Q_\alpha$, and contracts to $Q_\alpha$ in $D$. Note that if $Q_\alpha \subsetneq  Q_\beta$, then $N_\beta$ contains $P'_\alpha$ but $N_\alpha$ does not contain $P'_\beta$. Thus $N_\alpha$ survives in $(R:P'_\beta)$, but both $N_\alpha$ and $N_\beta$ blow up in $(R:P'_\alpha)$ since  both are principal and properly contain $P'_\alpha$.

Since $R$ is Pr\"ufer,  $(R:P'_\alpha)$ is the intersection $R_{P'_\alpha}$ and the localizations of $R$ at the maximal ideals that do not contain $P'_\alpha$ \cite[Theorem 3.1.2]{FHPbook}. Since  $R\subsetneq S$ are Pr\"ufer domains with the same quotient field,  $R_{P\cap R}=S_P$ for each (nonzero) prime $P$ of $S$. In particular,  $R_{P'_\alpha}=S_{P_\alpha}$ and if $P_\beta$ is a maximal ideal of $S$, then $S_{P_\beta}$ contains $V_\beta$. Thus $(R:P'_\alpha)= (S:P_\alpha)\cap F(P'_\alpha)$ where $F(P'_\alpha)=\bigcap \{V_\gamma\mid N_\gamma+P'_\alpha=R\}$. 

(7) Note that $P'_\alpha=\bigcap N_\alpha^n$ which implies that $P'_\alpha$ is a divisorial prime ideal of $R$ since $N_\alpha$ is invertible in $R$. If $P_\alpha$ is not a maximal ideal of $S$, then each ideal between  $(P_\alpha)^v$ and $P_\alpha$ (properly containing $P_\alpha$) survives in $(S:P_\alpha)$ and properly contains $P_\alpha$. In particular, there is a prime $P$ of $(S:P_\alpha)$ that properly contains $P_\alpha$. The contraction of $P$ to $(R:P'_\alpha)$ is a prime that properly contains $P'_\alpha$. Hence $P'_\alpha$  is not antesharp as an ideal of $R$. \end{proof}

For examples of primes satisfying cases 8 and 11 of Theorem \ref{th:1.12}, we use the ring of integer-valued polynomials  $\textrm{Int}(\mathbb Z)=\{f\in  \mathbb Q[X] \mid f(n)\in \mathbb Z$ for each $n\in \mathbb Z\}$ as a starting point. What makes $\textrm{Int}(\mathbb Z)$ particularly useful for our examples is the fact that it is a completely integrally closed two-dimensional Pr\"ufer domain  \cite[Sections V.2 and VI.2]{CCbook}.  Thus for each nonzero prime $P$, $(P:P)=\textrm{Int}(\mathbb Z)$ and $P$ is  branched. In fact, it is known that each  nonzero prime of $\textrm{Int}(\mathbb Z)$ has a trivial dual  \cite[Section VIII.5]{CCbook}. Also, no nonzero prime is idempotent. Thus by Corollary \ref{cor:1.4} and Proposition \ref{pr:1.5}, no nonzero prime  ideal is sharp and only the maximal ones are antesharp.  Each height-one nonmaximal prime is of the form $\mathfrak{P}_f:=f\mathbb Q[X]\cap D$ where $f\in \mathbb Q[X]$ is an irreducible polynomial (see \cite[Proposition V.2.7]{CCbook}). For our purposes we use $f(X)=X$.

\begin{example} $\boldsymbol{\Lambda}(P) =
\langle \nots,\nota, \notd,{\bf B},\noti\rangle$ \ (Case 11  of Theorem \ref{th:1.12}). 

Let   $D:= \textrm{Int}(\mathbb Z)$ and  consider the prime  ideal $P:=\mathfrak{P}_X$.  Since   $D_P= \mathbb Q[X]_{(X)}$, $PD_P$ is principal. Thus $P$ is branched and not idempotent. It is not maximal, since (for example) it is properly contained in the maximal ideal $M =\{f\in \textrm{Int}(\mathbb Z) \mid f(0)\in 2\mathbb Z\}$. Thus by the remarks above, $\boldsymbol{\Lambda}(P)=\langle \nots,\nota,\notd,{\bf B},\noti\rangle$.
 \end{example}

\begin{example} $\boldsymbol{\Lambda}(P) =
\langle \nots,\nota,{\bf D}, {\bf B},\noti\rangle$ \ (Case 8  of Theorem \ref{th:1.12}). 

Start with the ring of integer-valued polynomials  $D:= \textrm{Int}(\mathbb Z)$ and the nonmaximal prime  $Q: = {\mathfrak{P}_X}$. As we remarked in the previous example, (with the new notation) we have $D_Q =\mathbb Q[X]_{(X)}$ and $Q^{-1} =D$. As in Theorem \ref{th:2.12}, let $S:=D(Y)$, $P:= Q(Y)$ and $W:=  \mathbb Q[Y]_{(Y)}$ (here using $Y$ in place of the ``$X$" in [\ref{th:2.12}]). In this situation $S_P = \mathbb Q(Y)[X]_{(X)}$ and $V: =\psi_P^{-1}(W) =\mathbb Q[Y]_{(Y)}+X\mathbb Q(Y)[X]_{(X)}$. Set  $R:=  V \cap S$. Then (by Theorem \ref{th:2.12}) $R$ is a B\'ezout domain and $P$ is a nonmaximal divisorial prime of $R$ such that $(P:P) =(R : P)=S$. Since $P$ is not maximal in $S$,  $P$ is not antesharp in $R$. Also it is clear that $P$ is branched and not idempotent in $R$ (and in $S$). Thus $\boldsymbol{\Lambda}_R(P)=\langle \nots,\nota,{\bf D},{\bf B},\noti\rangle$.
 \end{example}
 
 The last four examples are obtained using the ring of entire functions $E$ as a base. The best single reference for the results needed for constructing the examples is M. Henriksen's paper \cite{Henprim} (other useful references are \cite{Henid}, \cite[pages 146-148 and Exercise 19 page 254]{Gilmit} and \cite[Section 8.1]{FHPbook}). We will provide some, but not all the details of the properties we need. Like the ring of integer-valued polynomials, $E$ is completely integrally closed. Unlike the ring of integer-valued polynomials, $E$ is a B\'ezout domain and it does have divisorial prime ideals, but each of these is a height-one principal maximal ideal of the form $M_a=(z-a)E$ for some uniquely determined $a\in \mathbb C$. Moreover,  all other nonzero primes have infinite height, some branched and some unbranched, each of these primes is contained in a unique maximal ideal and only the maximal ones of these are not idempotent. Since $E$ is a completely integrally closed Pr\"ufer domain, $P^{-1}=(P:P)=E$ except when $P$ is an invertible, necessarily height-one maximal, ideal (in which case $P^{-1}\supsetneq (P:P)=E$). Also, $E=\bigcap \{E_{M_a}\mid a\in \mathbb C\}$. Hence the height-one maximal ideals are the only sharp primes. 

The examples will be based on two specific types of nonmaximal primes of $E$.  First a little notation. For a nonzero  $g\in E$, we let $Z(g):=\{a\in \mathbb C\mid g(a)=0\}$. The set $Z(g)$ is nonempty if and only if $g$ is not a unit. It is also the case that $Z(g)$ is countable with no limit points. For each $a\in \mathbb C$,  we let $\mu_g(a)$ denote the multiplicity of $a$ as a zero of $g$ (equal to 0 if $g(a)\ne 0$, equal to $n\ (\ge 1)$ if $g(z)/(z-a)^n$ is in $E$ but  $g(z)/(z-a)^{n+1}$ is not). Two entire functions $g$ and $h$ are associates if and only if both $Z(g)=Z(h)$ and $\mu_g(a)=\mu_h(a)$ for each $a$. For convenience we set $\mu_r(a)=\infty$ if $r$ is the zero function.  
By a theorem of Weierstrass, there  are entire functions $e$ and $f$  such that $Z(e)=Z(f)$ is the set of positive integers with $\mu_e(n)=1$ and $\mu_f(n)=n$ for each positive integer $n$. Moreover, for any nonempty subset $A$ of $Z(e)$ and any sequence of nonnegative integers $\{n_k\}$, there is an entire function $g$ such that $Z(g)=A$ with $\mu_g(k)=n_k$ for each positive integer $k$ (with $n_k=0$ when $k\notin A$). Let $\mathcal U$ be a nonprincipal ultrafilter on the positive integers (so no subset in $\mathcal U$ is finite). For $r,s\in E$, set $r\equiv s$ if for each pair of positive integers $n$ and $m$, there are sets $U_n$ and $U_m$ in $\mathcal U$ such that $n\mu_r(k)\ge \mu_s(k)$ for each $k\in U_n$ and $m\mu_s(j)\ge \mu_r(j)$ for each $j\in U_m$. Based on this equivalence relation, the resulting equivalence classes are totally ordered using $[r]>[s]$ if for each $n$, there is a $U_n\in \mathcal U$ such that $\mu_r(k)>n  \mu_s(k)$ for each $k\in U_n$.   For $r\in E$ with $[r]>[e]$, there  are pairs of entire functions $a$ and $c$ such that $r$ and $a^2c$ are associates with $[r]=[a]$ (for example, there are entire functions $a$ and $c$ such that  $\mu_a(z)=\lfloor \mu_r(z)/2\rfloor$ for each $z$ and $\mu_c(z)=1$ whenever $\mu_r(z)$ is odd, then $r$ and $a^2c$  are associates with $[r]=[a]$).

The set $M:=\{ r\in E\mid [r]\ge [e]\}$ is a maximal ideal of $E$ that is not idempotent (in particular, $e$ is not in $M^2$) and it is minimal over the principal ideal $(e)$. Also  the set $P(M):=\{r\in E\mid [r]>[e]\}$ is a nonmaximal prime equal to $\bigcap_{n\ge 1} M^n$.  Similarly, $\overline P_f:=\{s\in E\mid [s]\ge [f]\}$ and $P_f:=\{s\in E\mid [s]>[f]\}$ are (nonmaximal) prime ideals with $M\supsetneq P(M)= P_e\supsetneq \overline P_f\supsetneq P_f=\bigcap_{n\ge 1} \overline P_f^n \supsetneq (0)$.   All four of $M$, $P(M)$, $\overline P_f$ and $P_f$ have infinite height (in fact, there is an uncountable descending chain of primes below $P_f$, and an uncountable chain between $P(M)$ and $\overline P_f$). Of these four, only $M$ is not idempotent (since $[r]>[e]$ implies the existence of a pair $a, c\in E$  such that $r=a^2c$ and $[r]=[a]$), both $M$ and $\overline P_f$ are branched ($\overline P_f$ is minimal over $(f)$) while both $P(M)$ and $P_f$ are unbranched. The latter follows from the fact that for any two nonzero entire functions $g$ and $h$ with $[g]>[h]\ge [e]$, there are nonzero  entire functions $b$ and $d$ such that $[b]>[g]>[d]>[h]$ (for example, let   $b,d\in E$ be such that $Z(b)=Z(g)\cap Z(e)$  with $\mu_b(k)=\mu_g(k)^2$ for each $k\in Z(b)$ and $Z(d)=Z(g)\cap Z(h)\cap Z(e)$ with $\mu_d(j)=\lfloor \sqrt {\mu_g(j)\mu_h(j)}\rfloor$ for each $j\in Z(d)$).   Thus $P(M)=\bigcup_{[r]>[e]} P_r$ and $P_f=\bigcup_{[g]>[f]}  P_g$. Also note that for $P\in \{M, P(M),\overline P_f, P_f\}$, $P^{-1}=(P:P)=E$. Thus only $M$ is antesharp, and as mentioned above, none are sharp.  
 
 As in the previous examples, the specific prime  of a certain type will be denoted ``$P$" in each of the examples below with $P$ built from either  $P(M)=P_e$ or $\overline P_f$, as defined above.

\begin{example} \label{ex:2.16} $\boldsymbol{\Lambda}(P) = \langle \nots,\nota,\notd, \notb,{\bf I}\rangle$ \ (Case 12  of Theorem \ref{th:1.12}). 

For this example, use the prime $P:=P(M)$ defined above. Then $P=P_e=\bigcup \{P_r\mid [r]>[e]\}$, so $P$ is unbranched, and therefore, idempotent. It is not maximal and  $P^{-1}=(P:P)=E$. Thus $P$ is not divisorial, not  antesharp and not sharp, so $\boldsymbol{\Lambda}(P)=\langle \nots,\nota,\notd,\notb,{\bf I}\rangle$. 
 \end{example}

\begin{example} $\boldsymbol{\Lambda}(P) = \langle \nots,\nota,{\bf D},\notb,{\bf I} \rangle$  (Case 9  of Theorem \ref{th:1.12}). 

For this example, set $Q:= P(M)=P_e$. Then  $\boldsymbol{\Lambda}_E(Q)=\langle \nots,\nota,\notd,\notb,{\bf I}\rangle$ as we remarked above. Thus by Theorem \ref{th:2.1}, the  nonmaximal prime ideal  $P:=Q(X)$ of $S:=E(X)$ has the same type as $Q$.  Let $V$ be the pullback of $W:=K[X]_{(X)}$ along $PS_P$ where $K:= E_Q/QE_Q$ and let  $R:=V\cap S$. Then (by Theorem \ref{th:2.12})  $P$ is a divisorial unbranched nonmaximal prime ideal of $R$ that is not antesharp. Hence $\boldsymbol{\Lambda}_R(P)=\langle \nots,\nota,{\bf D},\notb,{\bf I}\rangle$. 
\end{example}

\begin{example} \label{ex:2.18} $\boldsymbol{\Lambda}(P) = 
\langle \nots,\nota,\notd,{\bf B},{\bf I} \rangle$   (Case 10  of Theorem \ref{th:1.12}). 

 In this example, let  $P:= \overline P_f=\{r\in E\mid [r]\ge [f]\}$ where $[f]>[e]$ with $f$ and $e$ as above. Then as we discussed above, $P$ is a nonzero nonmaximal prime that is branched,  idempotent and  not divisorial, so $\boldsymbol{\Lambda}(P)=\langle \nots,\nota,\notd,{\bf B},{\bf I}\rangle$. 
\end{example} 

We finish the examples with a nonzero nonmaximal prime $P$ that fits case 7 of Theorem \ref{th:1.12}, $\boldsymbol{\Lambda}(P)=\langle \nots,\nota,{\bf D},{\bf B},{\bf I}\rangle$.

\begin{example} $\boldsymbol{\Lambda}(P) = 
\langle \nots,\nota,{\bf D},{\bf B},{\bf I} \rangle$  \ (Case 7  of Theorem \ref{th:1.12}).

 As in the previous example, start with  $Q:=\overline P_f$. In $S:=E(X)$, the prime $P= QS$ has the same type in $S$ as $Q$ has in $E$, namely $\boldsymbol{\Lambda}_S(P)=\langle \nots,\nota,\notd,{\bf B},{\bf I}\rangle$.  Let $W$  and $V$ be as in Theorem \ref{th:2.12}. Then  in $R:=V\cap S$, the prime (nonmaximal) ideal $P$ is still branched and idempotent and neither sharp nor antesharp, but it is divisorial in $R$. Hence $\boldsymbol{\Lambda}_R(P)=\langle \nots,\nota,{\bf D},{\bf B},{\bf I}\rangle$.
\end{example}

\end{document}